\documentclass[journal]{IEEEtran}
\usepackage[pdftex,dvipsnames]{xcolor}
\usepackage{tikz}
\usepackage{tikz}
\usepackage{array}
\newcolumntype{P}[1]{>{\centering\arraybackslash}p{#1}} 
\ifCLASSOPTIONcompsoc
\usepackage[caption=false,font=normalsize,labelfont=sf,textfont=sf]{subfig}
\else
\usepackage[caption=false,font=footnotesize]{subfig}
\fi
\usepackage{textcomp}
\usepackage{stfloats}
\usepackage{verbatim}
\usepackage{dirtytalk}
\usepackage{multirow}
\usepackage{rotating}
\usepackage{booktabs}
\usepackage{tikz}
\usepackage{graphicx}
\usepackage{amsmath,amsfonts}
\usepackage{bm}
\usepackage{hyperref}
\usepackage{url}
\usepackage{cite}
\usepackage{subfig}
\usepackage{algorithm} 
\usepackage{algpseudocode} 
\usepackage{tikz-network}
\usetikzlibrary{arrows,automata,positioning}

\usepackage{tikz}
\usetikzlibrary{trees}
\usetikzlibrary{chains, fit}
\usepackage{listings}
\usepackage{xcolor,colortbl}
\usetikzlibrary{shapes,arrows}
\usetikzlibrary{shadings}
\usetikzlibrary{automata,topaths}
\usetikzlibrary{shapes}
\usetikzlibrary{shapes.geometric, arrows}
\usepackage{tikz-3dplot}
\usetikzlibrary{matrix}
\usepackage{comment}

\usepackage{pgfplots}
\usetikzlibrary{pgfplots.groupplots}

\newcommand{\R}{{\rm I\!R}}

\newcommand{\addb}[1]{{{\color{blue!0!black}#1}}}

\newcommand{\addw}[1]{{{\color{white!100!black}#1}}}
\usetikzlibrary{quotes}

\ifCLASSINFOpdf
\else
\fi
\newcommand{\boundellipse}[3]
{(#1) ellipse (#2 and #3)
}

\newcommand{\cf}{{\emph{cf.~}}}
\newcommand{\tran}{^{\mbox{\scriptsize T}}}
\usepackage{pifont}
\newtheorem{theorem}{Theorem}
\newtheorem{lemma}[theorem]{Lemma}
\newtheorem{assump}{Assumption}
\newtheorem{assump2}{Assumption}
\numberwithin{assump2}{section} 

\newtheorem{remark}{Remark}

\newtheorem{prop}{Proposition}
\newtheorem{prop2}{Proposition}
\numberwithin{prop2}{section}

\newtheorem{corollary}{Corollary}

\makeatletter
\newenvironment{breakablealgorithm}
  {
   \begin{center}
     \refstepcounter{algorithm}
     \hrule height.8pt depth0pt \kern2pt
     \renewcommand{\caption}[2][\relax]{
       {\raggedright\textbf{\fname@algorithm~\thealgorithm} ##2\par}%
       \ifx\relax##1\relax 
         \addcontentsline{loa}{algorithm}{\protect\numberline{\thealgorithm}##2}%
       \else 
         \addcontentsline{loa}{algorithm}{\protect\numberline{\thealgorithm}##1}%
       \fi
       \kern2pt\hrule\kern2pt
     }
  }{
     \kern2pt\hrule\relax
   \end{center}
  }
\makeatother
\makeatletter
\newcommand{\nosemic}{\renewcommand{\@endalgocfline}{\relax}}
\newcommand{\dosemic}{\renewcommand{\@endalgocfline}{\algocf@endline}}
\let\oldnl\nl
\newcommand{\nonl}{\renewcommand{\nl}{\let\nl\oldnl}}
\makeatother
\usepackage{amssymb}
\newcommand{\QE}{\hfill $\blacksquare$}

\begin{document}
\bstctlcite{IEEEexample:BSTcontrol}
\title{On the Primal Feasibility in Dual Decomposition Methods Under Additive and Bounded Errors}
\author{Hansi~Abeynanda,
        Chathuranga~Weeraddana,
        and~G. H. Jayantha~Lanel
}



\author{Hansi~Abeynanda,
        Chathuranga~Weeraddana,
        G. H. J.~Lanel,
        and~Carlo~Fischione
\thanks{H. Abeynanda is with the Mathematics Unit of Sri Lanka Institute of Information Technology, Sri Lanka (e-mail: kavindika.a@sliit.lk).}
\thanks{C. Weeraddana is with Centre for Wireless Communication,
University of Oulu, Finland (e-mail: chathuranga.weeraddana@oulu.fi).}
\thanks{G. H. J. Lanel is with the Department of Mathematics, University of Sri Jayewardenepura, Sri Lanka (e-mail: ghjlanel@sjp.ac.lk).}
\thanks{C. Fischione is with the School of Electrical Engineering and Computer Science, KTH Royal Institute of Technology, Sweden (e-mail: carlofi@kth.se).}
}
\maketitle
\begin{abstract}
With the unprecedented growth of signal processing and machine learning application domains, there has been a tremendous expansion of interest in distributed optimization methods to cope with the underlying large-scale problems. Nonetheless, inevitable system-specific challenges such as limited computational power, limited communication, latency requirements, measurement errors, and noises in wireless channels impose restrictions on the exactness of the underlying algorithms. Such restrictions have appealed to the exploration of algorithms’ convergence behaviors under inexact settings. Despite the extensive research conducted in the area, it seems that the analysis of convergences of dual decomposition methods concerning primal optimality violations, together with dual optimality violations is less investigated. Here, we provide a systematic exposition of the convergence of \emph{feasible points} in dual decomposition methods under inexact settings, for an important class of global consensus optimization problems. Convergences and the rate of convergences of the algorithms are mathematically substantiated, not only from a dual-domain standpoint but also from a primal-domain standpoint. Analytical results show that the algorithms converge to a neighborhood of optimality, the size of which depends on the level of underlying distortions. 
\end{abstract}
%
\IEEEpeerreviewmaketitle



\section{Introduction}\label{sec:Introduction}


%
A prevalent problem in many 
application fields, including signal processing, machine learning, telecommunication networks, control systems, robotics, and network applications, among others, \cite{Palomar-book-2010,Carlo-Machine-Learning,Nedic-conv-rate-book-2015,Yang-survey-of-distributed-optimization-2019,Nedic-Distributed-Optimization-for-Control-2018,Boyd-resource-allocation-dual-decomposition-2004,Madan-2006,Nedic-2009,Zhang-2021,Daniel-Survey-Distributed-Optimization-2017} is
 \begin{equation} \label{eq:main-problem}
\begin{array}{ll}
\mbox{minimize} & \sum_{i{=}1}^{m}f_i(\mathbf{z}) \\
\mbox{subject to} & \mathbf{z}\in\mathcal{Y}, 
\end{array}
\end{equation}
where the variable is $\mathbf{z}\in\R^n$ and $\mathcal{Y}\subseteq\R^{n}$ is considered as a convex and closed constraint set. Each $f_i:\R^n\to \R$ is a strictly convex and closed function associated with subsystem $i$. Here $\mathbf{z}$ is called the \emph{public} variable. The problem is sometimes known as the consensus problem. Other real-world applications of~\eqref{eq:main-problem} include networked vehicles, smart power grids, control of UAV/multiple robots, and TCP control systems, \cite{halsted-Multiple-robots-2021,Javier-multiple-robot-2016}. 
In practice, the unprecedented growth of the size of modern datasets, decentralized collection of datasets, and underlying high-dimensional decision spaces, prevents the applicability of centralized methods such as interior-point algorithms~\cite{Boyd-Interior-point-methods} for solving~\eqref{eq:main-problem}. They entail the development of salable distributed algorithms~\cite{Boyd-Parikh-Chu-Peleato-Eckstein-2010,Boyd-Parikh-Proximal-Algorithms-2014}.


Two commonly used first-order algorithms are dual decomposition methods~\cite{Low-Layering-as-Optimization-Decomposition-2007,Palomar-2007,R2-Nedic-Dual-Approach-Optimal-Algorithms-2020,R3-Scaman-2017,R4-Hendrikx-2020} and approaches coalescing consensus algorithms with subgradient methods \cite{Yang-survey-of-distributed-optimization-2019,Nedic-conv-rate-book-2015}. 
The \emph{Simplicity} of the implementation of dual decomposition techniques qualifies it as a promising solution method, especially in large-scale optimization problems~\cite[\S~10]{Palomar-book-2010}\cite[\S~3.4]{Bertsekas-book-1997}. In almost all distributed methods, underlying applications consist of subsystems making local decisions and coordinating information to solve problem~\eqref{eq:main-problem}. During this process, systems have to operate under many nonideal settings. Thus, the analysis of algorithms with nonideal settings has been an appealing area of study
\cite{Polyak-Intro-Opt-1987,Nedic-quantization-effects-2008,Peng-Quantized-subgradient-algorithm-2014,Chang-quantization-grid-shrink-2018,Solodov-1998,Rabbat-Quantized-incremental-algorithms-2005,Nedic-2010,Nestero-First-order-methods-2014,Bertsekas-1999,Chen-biased-SGD-2018,Ahmad-biased-SGD-2021,Sindri-2021,Sindri-2018,Sindri-2020,R1-Liu-2021,Necoara-inexact-dual-decomposition-2014,Yifan-inexact-dual-decomposition-2021,sindri_Communication_Complexity_2018}.

The main focus of the manuscript resides in \emph{dual decomposition} with inexact gradient methods for~\eqref{eq:main-problem} with a greater emphasis on convergences from a \emph{primal-domain} standpoint. A wide range of distortions is considered where nothing except their norm boundedness is assumed. We explicitly provide the convergences and the rate of convergences of related algorithms concerning the dual optimality violations and the primal optimality violations. More importantly, we establish the convergences of primal feasible points with a complete convergence rate analysis.

\subsection{Related Work}

An elegant discussion on the influence of noise in subgradient methods can be found in \cite{Polyak-Intro-Opt-1987} under both differentiable and nondifferentiable settings (e.g., \cite[\S~4, and \S~5.5]{Polyak-Intro-Opt-1987}). More importantly, \cite{Polyak-Intro-Opt-1987} provides a repertory of techniques that can serve as building blocks that are indispensable when analyzing algorithms with imperfections.

Algorithms based on combining consensus algorithms with subgradient methods under nonideal settings have been discussed in \cite{Nedic-quantization-effects-2008,Peng-Quantized-subgradient-algorithm-2014,Chang-quantization-grid-shrink-2018}. Under assumptions such as uniform boundedness of subgradients, the convergences of underlying algorithms are derived in \cite{Nedic-quantization-effects-2008} and \cite{Peng-Quantized-subgradient-algorithm-2014}. However, the boundedness assumption might restrict the range of applicability of the methods. For example, in many applied fields, it is now commonplace to form the objective function with a quadratic regularization term, where the bounded assumption is no longer affirmative. Moreover, in~\cite{Peng-Quantized-subgradient-algorithm-2014,Chang-quantization-grid-shrink-2018} authors consider distortions due to quantization, which are diminishing. Although a diminishing error is favorable from a standpoint of establishing desirable convergences, it
cannot capture distortions that are persistent, e.g., measurement errors.



Inexact gradient/subgradient methods are considered in \cite{Solodov-1998,Bertsekas-1999,Rabbat-Quantized-incremental-algorithms-2005,Nedic-2010,Nestero-First-order-methods-2014,Chen-biased-SGD-2018,Ahmad-biased-SGD-2021,Sindri-2021}. It is usually the case that subgradient type algorithms are used to solve dual problems in a dual decomposition setting. The effect of noise in subgradient type methods has been discussed in \cite{Solodov-1998,Rabbat-Quantized-incremental-algorithms-2005,Nedic-2010,Nestero-First-order-methods-2014} with \emph{compact} constraint sets. From a distributed optimization standpoint with dual decomposition, compactness is a restriction because constraint sets appearing in the dual-domain usually turn out to be noncompact. The errors considered in \cite{Bertsekas-1999} are diminishing, which is a more restrictive property, as we have already pointed out. In a machine learning setting, \cite{Chen-biased-SGD-2018} consider distortions that are regulated by an adequate choice of the sample size. However, the method is not amenable in the dual decomposition setting when the underlying errors are not necessarily controllable. References \cite{Ahmad-biased-SGD-2021,Sindri-2021} appear to be readily applied in a distributed optimization setting with dual decomposition.

Inexact gradient methods~\cite{Sindri-2018,Sindri-2020,Yifan-inexact-dual-decomposition-2021,R1-Liu-2021,Necoara-inexact-dual-decomposition-2014,sindri_Communication_Complexity_2018} considered within the dual-domain are closely related to our study. A concise summary of specific technical assumptions and related convergence rate results, among others are tabulated in Table~\ref{Table:comparison-table}. Broadly speaking, the authors in \cite{Sindri-2018} use normalized and quantized gradients, which admits a \emph{zooming-in and quantized} policy. A related \emph{zooming-in and quantized} policy is discussed also in \cite{Sindri-2020}. The modeling assumptions in such a policy are restrictive to be adopted in general since it imposes conditions on distortions to diminish as the iteration number increases. The inexactnesses of dual gradients considered in \cite{R1-Liu-2021} are due to the accuracy of subproblem solvers in an attempt to save communication during subsystem coordination. Technically, the underlying distortions are again diminishing, as in \cite{Sindri-2018,Sindri-2020}. It is worth noting that the references \cite{Sindri-2018,Sindri-2020,R1-Liu-2021} have not considered any primal-domain convergences.

In references~\cite{Yifan-inexact-dual-decomposition-2021,Necoara-inexact-dual-decomposition-2014,sindri_Communication_Complexity_2018}, specific results associated with primal-domain convergences have been established within a dual decomposition setting. Authors in~\cite{Necoara-inexact-dual-decomposition-2014} establish convergence results in the primal-domain with strong convexity assumptions on objective functions and certain norm boundedness assumptions on the constraint functions. It is worth emphasizing that the sequence of primal variables pertaining to the convergence results in ~\cite{Necoara-inexact-dual-decomposition-2014} does not admit primal feasibility, in general. Thus, the sequence is indeed \emph{infeasible}, despite the number of iterates the algorithm is performed. However, their results characterize an asymptotic \emph{primal optimality violation}, i.e., a sequence converging to a ball around the optimality at a rate of $O(1/k)$. 
Similar convergence results of {primal optimality violation} have been derived in \cite{Yifan-inexact-dual-decomposition-2021,sindri_Communication_Complexity_2018}, where only the problems of the form of sharing are considered. Reference \cite{Yifan-inexact-dual-decomposition-2021} establishes convergence to a neighborhood of optimality at a rate of $O(1/\sqrt{k})$ with strong convexity and twice-differentiability assumptions on objective functions and compactness assumptions on the constraint set. In contrast, a linear convergence rate is established in \cite{sindri_Communication_Complexity_2018} with more restricted assumptions. For example, the associated dual gradients of the considered problem in~\cite{sindri_Communication_Complexity_2018} are scalars, and when quantizing, additional assumptions on the compactness of the dual-domain are artificially imposed for tractability. Moreover, diminishing distortions are considered due to quantization. The authors in \cite{sindri_Communication_Complexity_2018} further provide a sequence of \emph{primal feasible points} that converges to the optimality at a rate of~$O(1/\sqrt{k})$.



\begin{table*}[t]
\caption{Comparison of the Results in This Paper with the Most Related Existing Results. Acronyms: DO $=$ Dual Optimality, PO $=$ Primal Optimality, PF $=$ Primal Feasibility, CR=Convergence Rate, DLAG $=$ Dual Accelerated Method with Lazy Approximate Gradient} \label{Table:comparison-table}
\begin{minipage}{\textwidth}
\begin{center}
\begin{tabular}{|m{5.3em}|m{2em}|m{9em}|m{4.5em}|c|c|c|c|c|}
 \hline
 \multirow{2}{*}{\hspace{0.25cm}Method} & \multirow{2}{*}{Ref.} & \multicolumn{2}{c|}{Main Assumptions} & \multicolumn{2}{c|}{Stepsize} & \multirow{2}{*}{DO/CR} & \multirow{2}{*}{PO/CR} & \multirow{2}{*}{PF/CR}\\ 
 \cline{3-6}
 &  & Primal functions & Errors & Constant & Nonsummable  & & &\\
 \hline
   \multirow{3}{5em}{Gradient descent centralized} & \multirow{3}{*}{\cite{Sindri-2018}} & \multirow{3}{9em}{Assump.~\ref{Assumption:Lipschitz-Convex-Grad-f}} &  \multirow{3}{*}{Diminishing} & YES & YES & YES &  NO & NO \\
 \cline{5-9}
 & &  &  &  \multicolumn{2}{c|}{Constant} & $O(1/\sqrt{k})$~\footnote{\label{note1}CR is given with respect to the norm of the gradient.}  &  NA & NA \\
 \cline{5-9}
 & &  &  &  \multicolumn{2}{c|}{Nonsummable} & CR not derived
 & NA & NA \\
 \hline
 \multirow{3}{6em}{Dual decomposition with gradient descent} & \multirow{2}{*}{\cite{Sindri-2020}} & \multirow{3}{9em}{Assump.~\ref{Assumption:Strongly-Convex-f} \& Assump.~\ref{Assumption:Lipschitz-Convex-Grad-f}} &  \multirow{3}{*}{Diminishing} & YES & NO & YES & NO & NO \\
 \cline{5-9}
 & &  & &  \multicolumn{2}{c|}{Constant} & Linear rate~\footnote{\label{note2}CR is given with respect to decision variables.} & NA  & NA \\
 \cline{5-9}
 & &  & &  \multicolumn{2}{c|}{Nonsummable} & NA & NA  & NA  \\
 \hline
 \multirow{3}{7em}{DLAG and Multi-DLAG} & \multirow{3}{*}{\cite{R1-Liu-2021}} & \multirow{3}{9em}{Assump.~\ref{Assumption:Strongly-Convex-f} \& Assump.~\ref{Assumption:Lipschitz-Convex-Grad-f}} &  \multirow{3}{*}{Diminishing} & YES  & NO  & YES & NO & NO \\
 \cline{5-9}
 & &  & & \multicolumn{2}{c|}{Constant} & Linear rate~\footnote{\label{note3}CR is given with respect to function values.} & NA  & NA \\
  \cline{5-9}
 & &  & & \multicolumn{2}{c|}{Nonsummable} & NA & NA  & NA \\
 \hline
 \multirow{3}{5em}{Gradient descent centralized} & \multirow{3}{*}{\cite{Necoara-inexact-dual-decomposition-2014}} & \multirow{3}{9.5em}{Assump.~\ref{Assumption:Strongly-Convex-f}, $f_i$s twice differentiable \& the Jacobian of the constraint function is bounded} & \multirow{3}{*}{Bounded} & NO  & YES  & YES & YES & YES\footnote{Primal feasibility is discussed for a very
special case of a model predictive control problem concerning the convergences of primal function values.} \\
 \cline{5-9}
 & &  & &  \multicolumn{2}{c|}{Constant} & NA & NA & NA \\
 \cline{5-9}
 & &  & &  \multicolumn{2}{c|}{Nonsummable} & $O(1/k)$~\textsuperscript{\ref{note3}}  & $O(1/k)$~\textsuperscript{\ref{note3}} & $O(1/\sqrt{k})$~\textsuperscript{\ref{note3}} \\
 [0.1cm]
 \hline
 \multirow{3}{6em}{Dual decomposition with gradient descent} & \multirow{3}{*}{\cite{Yifan-inexact-dual-decomposition-2021}} & \multirow{3}{8.7em}{Assump.~\ref{Assumption:Strongly-Convex-f}, $f_i$s twice differentiable, \& feasible set compact} & \multirow{3}{*}{Bounded} & YES & NO & YES & YES & NO \\
 \cline{5-9}
 & &  & &  \multicolumn{2}{c|}{Constant} & $O(1/\sqrt{k})$~\textsuperscript{\ref{note3}} & $O(1/\sqrt{k})$~\footnote{\label{note4}CR is given with respect to both decision variables and function values.} & NA  \\
 \cline{5-9}
 & &  & &  \multicolumn{2}{c|}{Nonsummable} & NA  & NA  & NA  \\
 \hline
  \multirow{3}{6em}{Dual decomposition with gradient descent} & \multirow{3}{*}{\cite{sindri_Communication_Complexity_2018}} & \multirow{3}{8em}{Assump.~\ref{Assumption:Strongly-Convex-f} \& dual-domain compact} & \multirow{2}{*}{Diminishing} & YES & NO & YES  & YES & YES\footnote{Convergence of primal feasible points is shown with more restricted assumptions.} \\
 \cline{5-9}
 & &  & &  \multicolumn{2}{c|}{Constant} & Linear rate~\textsuperscript{\ref{note3}} & Linear rate~\textsuperscript{\ref{note2}}  & $O(1/\sqrt{k})$~\textsuperscript{\ref{note2}} \\
 \cline{5-9}
 & &  & &  \multicolumn{2}{c|}{Nonsummable} & NA & NA  & NA \\
 \hline
 \multirow{6}{6em}{Dual decomposition with gradient descent} & \multirow{6}{2em}{This paper} & \multirow{3}{*}{Assump.~\ref{Assumption:Strongly-Convex-f}} &  \multirow{3}{*}{Bounded} & YES & YES & YES & YES & YES \\
 \cline{5-9}
 & & & & \multicolumn{2}{c|}{Constant} & $O(1/\sqrt{k})$~\textsuperscript{\ref{note1}}& $O(1/\sqrt[4]{k})$~\textsuperscript{\ref{note4}} &  $O(1/\sqrt[4]{k})$~\textsuperscript{\ref{note2}}\\
 \cline{5-9}
 & & & & \multicolumn{2}{c|}{Nonsummable} &
 $O(1/\sqrt{k^{1-p}})$~\textsuperscript{\ref{note1}} & $O(1/\sqrt[4]{k^{1-p}})$~\textsuperscript{\ref{note4}} & $O(1/\sqrt[4]{k^{1-p}})$~\textsuperscript{\ref{note2}}\\
 \cline{3-9}
 & & \multirow{3}{9em}{Assump.~\ref{Assumption:Strongly-Convex-f} \& Assump.~\ref{Assumption:Lipschitz-Convex-Grad-f}} &  \multirow{3}{*}{Bounded} & YES & YES & YES & YES & YES \\
 \cline{5-9}
 & &  & &  \multicolumn{2}{c|}{Constant} & Linear rate~\textsuperscript{\ref{note3}} & Linear rate~\textsuperscript{\ref{note4}} & Linear rate~\textsuperscript{\ref{note4}} \\
 \cline{5-9}
 & &  & &  \multicolumn{2}{c|}{Nonsummable} & $O(1/k^{c/p})$~\textsuperscript{\ref{note3}} & $O(1/k^{c/2p})$~\textsuperscript{\ref{note4}} & $O(1/k^{c/2p})$~\textsuperscript{\ref{note4}} \\
 \hline
\end{tabular}
\end{center}
\end{minipage}
\label{Table:summary-table-Case1} 
\end{table*}

Under inexact settings, the literature relying on dual decomposition methods predominantly considers convergences in dual-domain. In contrast, primal-domain convergences have received relatively little attention. Even though there are a few works establishing convergence results from a primal optimality violation standpoint still under restricted assumptions, no primal feasibility
guarantees of such sequences are established. Convergence results of \emph{primal feasible} sequences have been less investigated though such convergence results are of vital importance in a multitude of practical application domains, see \cite{Sindri-feasible-methods-2021} and references therein. Therefore, it is desirable to have an exposition that lays out primal-domain convergence results of dual decomposition methods under inexact settings.


\subsection{Our Contribution}

A problem of minimizing a global function which is a sum of local convex objective functions, under general convex constraints is considered [\cf problem~\eqref{eq:main-problem}]. 
Within this setting, the main contributions of the paper are summarized below.
\begin{enumerate}
    \item \emph{Inexact distributed algorithms:} A fully distributed algorithm is proposed  based on dual decomposition techniques [\cf~\S~\ref{sec:Imperfect-Coordination}] and gradient methods. The proposed algorithm is inexact in the sense that the underlying subproblem coordination is imperfect. Our modeling captures a wide range of distortions, including quantization errors, approximation errors, errors due to subproblem solver accuracy, noise in wireless settings, and measurement errors, among others, as long as they are \emph{additive} and  \emph{bounded} [\cf~\S~\ref{sec:Imperfect-Coordination}, Remark~\ref{Remark:Genarality-of-Error}]. The algorithm is analyzed  under \emph{fixed} and \emph{nonsummable} stepsize rules. 
    \item \emph{Primal optimality violation:} Under mild conditions, the convergence of the algorithm in the primal-domain is analytically substantiated [\cf~\S~\ref{sec:Convergence-Analysis-I}, Proposition~\ref{Proposition:Lipschitz-continuous-g-Primal-Result} and Proposition~\ref{Proposition:Strongly-convex-Lipschitz-continuous-constant-step-Primal}]. Despite \emph{primal infeasibility}, we show that the algorithms get into a neighborhood of optimality, the size of which depends on the level of underlying distortions. Convergence rates are also derived.
    \item \emph{Constructing primal feasible points and their optimality:} How to construct a feasible solution by using current infeasible primal variables is highlighted [\cf~\S~\ref{sec:Convergence-Analysis-I}, \cf~Remark~\ref{Remark:Compute-Feasible-Point}]. Under mild conditions, convergences of the algorithms in the primal-domain, while maintaining feasibility, are established [\cf~\S~\ref{sec:Convergence-Analysis-I}, Proposition~\ref{Proposition:Feasible-Point-Lipschitz-continuous-g-Primal-Result}, Proposition~\ref{Proposition:Feasible-Point-convex-Lipschitz-continuous-constant-step-Primal}, and Remark~\ref{Remark:Feasible-Point-convex-Lipschitz-continuous-nonsumable-step-Primal}]. Convergence rates are also derived and concisely compared with the most relevant state-of-the-art, \cf~Table~\ref{Table:comparison-table}.
\end{enumerate}

\subsection{Notation}


We use normal font lowercase letters $x$, bold font lowercase letters $\mathbf{x}$, bold font uppercase letters $\mathbf{X}$, and calligraphic font $\mathcal{X}$ to denote scalars, vectors, matrices, and sets, respectively.
For a given matrix $\mathbf{X}$, $\mathbf{X}\tran$ denotes the matrix transpose. Matrices $\mathbf{I}_{n}$ and $\mathbf{1}_{{m\times n}}$ denote the $n\times n$ identity matrix and the $m\times n$ matrix with all entries equal to one, respectively. The Kronecker product of matrices $\mathbf{A}$ and $\mathbf{B}$ is denoted by $\mathbf{A}\otimes \mathbf{B}$. 
The set of real numbers, set of extended real numbers, set of real $n$-vectors, set of real $m\times n$ matrices, set of positive integers, and set of nonnegative integers are denoted by $\R$, $\overline{\R}$, $\R^n$, $\R^{m\times n}$, $\mathbb{Z}_+$, and $\mathbb{Z}^0_+$, respectively. For $\mathbf{x}\in\R^n$, $\|\mathbf{x}\|$ and for $\mathbf{A}\in\R^{m\times n}$, $\|\mathbf{A}\|$ denote the $\ell_2$-norm and matrix 2-norm, respectively. The domain of a function $f:\R^n\to \R$ is a subset of $\R^n$ and is denoted by $\texttt{dom}~f$. The asymptotic notations \say{big Oh} and \say{small oh} are denoted by $O(\cdot)$ and $o(\cdot)$, respectively.

\subsection{Organization of the Paper}

The rest of the paper is organized as follows. In \S~\ref{sec:Dual-Decomposition-Applied-to-1}, we
apply dual decomposition techniques to our main problem~\eqref{eq:main-problem}.
The modeling of errors that underlie the inexactness, together with the proposed algorithm, is presented in \S~\ref{sec:Imperfect-Coordination}.
In \S~\ref{sec:Convergence-Analysis-I}, the main convergence results are presented. Numerical experiments are presented in \S~\ref{sec:Numerical-Results}. Lastly,
\S~\ref{Conclusion} concludes the paper, followed by appendices.

\section{Dual Decomposition Applied to~\eqref{eq:main-problem}}\label{sec:Dual-Decomposition-Applied-to-1}

\begin{figure}[t!]
\centering
\begin{tikzpicture}
\tikzset{vertex/.style = {shape=circle,draw,minimum size=0.2 cm}}
\tikzset{edge/.style = {->,> = latex}}
\node[vertex] (1) at  (0,0) {$1$};
\node[vertex] (2) at  (2,0) {$2$};
\node[vertex] (3) at  (4,0) {$3$};
\node[vertex] (4) at  (6,0) {$m$};
\node[vertex] (5) at  (3,-1.5) {$\mathbf{z}$};
\draw [thick, dotted]  (4.5,0) -- (5.5,0);
\draw [thick]  (1.330) -- node[left, xshift=-0.4cm, yshift=0.3em] {$\mathbf{y}_{1}$} (5.180);
\draw [thick]  (2.300) -- node[left, xshift=0em, yshift=0em] {$\mathbf{y}_{2}$} (5.120);
\draw [thick]  (3.240) -- node[right, xshift=0em, yshift=0em] {$\mathbf{y}_{3}$} (5.60);
\draw [thick]  (4.210) -- node[right, xshift=0.4cm, yshift=0.3em] {$\mathbf{y}_{m}$} (5.0);
\end{tikzpicture}
\caption{Decomposition Structure: There are $m$ SSs with the public variable $\mathbf{z}$. Functions associated with SSs are $f_i(\mathbf{z})$, $i\in\{1,\ldots,m\}$.}
\label{Fig:Decomposition-Structure}
\end{figure}
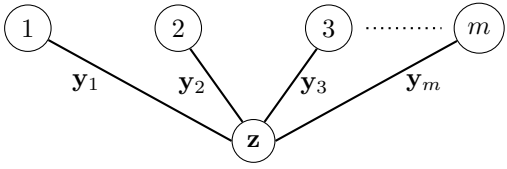

A commonly used technique to yield distributed solution methods for problem~\eqref{eq:main-problem} is based on dual decomposition~\cite{Boyd-EE364b-PrimDualDecomp-07}, where the decomposition structure of the underlying problem places a crucial role, see \figurename~\ref{Fig:Decomposition-Structure}.

We start by associating with each subsystem (SS)~$i$, a \emph{private variable} $\mathbf{y}_{i}$, together with necessary constraints to ensure their consistency $\mathbf{z}=\mathbf{y}_{i}$ for all $i=1,\ldots,m$. It is commonplace to impose the consistency of private variables as $\mathbf{y}_{i}=\mathbf{y}_{{i+1}}, \ i=1,\ldots,m-1$~\cite{Boyd-Parikh-Chu-Peleato-Eckstein-2010,Boyd-Parikh-Proximal-Algorithms-2014,Boyd-EE364b-PrimDualDecomp-07}. Thus, problem~\eqref{eq:main-problem} is equivalently reformulated as follows:
\begin{equation} \label{eq:distributed-problem}
\begin{array}{ll}
\mbox{minimize} & f(\mathbf{y})=\sum_{i{=}1}^{m}f_i(\mathbf{y}_{i}) \\
\mbox{subject to} & \mathbf{y}_{i}\in\mathcal{Y},\ i=1,\ldots,m\\ 
& \mathbf{y}_{i}=\mathbf{y}_{{i+1}}, \ i=1,\ldots,m-1,
\end{array}
\end{equation}
where $\mathbf{y}_{i}\in\R^n$, $i=1,\ldots,m$, are newly introduced local versions of the public variable $\mathbf{z}$ and $\mathbf{y}=[\mathbf{y}_{1}\tran \ \ldots \ \mathbf{y}_{m}\tran]\tran$.

In general, the consistency can be imposed by $\sqrt{\bar{\mathbf{W}}}\mathbf{y}=\mathbf{0}$ for some appropriately chosen matrix $\bar{\mathbf{W}}$.\footnote{For example, problem~\eqref{eq:main-problem} considered over a general communication graph admits $\bar{\mathbf{W}}=\mathbf{W}\otimes\mathbf{I}_n$, where $\mathbf{W}$ is the Laplacian of the graph \cite{R1-Liu-2021,Sindri-2020}.} It is worth pointing out that despite how the consistency of private variables is imposed, the essence of our technical derivations in the rest of the paper remains intact, \cf~Appendix~\ref{Appendix:General-communication-graph}. Hence, in the following, without loss of generality, we consider the consistency as specified in problem~\eqref{eq:distributed-problem} for notational convenience.

Note that the objective function of \eqref{eq:distributed-problem} is now separable. 
%
Let $\boldsymbol{\lambda}_{i}\in\R^n$ denote the Lagrange multiplier associated with the constraint $\mathbf{y}_{i}=\mathbf{y}_{{i+1}}$, $i=1,\ldots,m-1$ and $\boldsymbol{\lambda}=[\boldsymbol{\lambda}_1\tran \ \ldots \ \boldsymbol{\lambda}_{m-1}\tran]\tran$ for clarity.  
Then, the dual function $g:\R^{n(m-1)}\rightarrow \overline{\R}$ corresponding to \eqref{eq:distributed-problem} is given by 

\begin{align} \label{eq:dual-function}
    g(\boldsymbol{\lambda)}& = {\underset{\mathbf{y}_{i}\in\mathcal{Y},\ i=1,\ldots,m}{\inf}}\left[\sum_{i=1}^{m}f_i(\mathbf{y}_{i}) {+} \sum_{i=1}^{m{-}1} \boldsymbol{\lambda}_i\tran(\mathbf{y}_{i}{-}\mathbf{y}_{{i+1}})\right] \allowdisplaybreaks\\ \label{eq:dual-function-split}
    & = \sum_{i=1}^{m} \ \underbrace{\underset{\mathbf{y}_{i}\in\mathcal{Y}}{\inf}~\left[f_i(\mathbf{y}_{i}){+}(\boldsymbol{\lambda}_i-\boldsymbol{\lambda}_{i-1})\tran \mathbf{y}_{i}\right].}_{\mbox{\small{subproblem}}~i}  \allowdisplaybreaks
\end{align}
The last equality follows because, for fixed $\boldsymbol{\lambda}$, the infimization can be performed in parallel by each SS~\footnote{Here we have $\boldsymbol{\lambda}_{0}=\boldsymbol{\lambda}_{m}=\mathbf{0}$.}. Thus, associated with each SS, there is a subproblem that can be handled locally. The dual problem is given by
\begin{equation}\label{eq:dual-problem}
    \underset{\boldsymbol{\lambda}\in\R^{n(m-1)}}{\text{maximize}} \quad g(\boldsymbol{\lambda}).
\end{equation}

\def\a{.4}

The dual decomposition furnishes a mechanism to coordinate subproblems, one for each SS [\cf~\eqref{eq:dual-function-split}],  to solve problem~\eqref{eq:dual-problem} by using an iterative algorithm. We note that the
dual function is always concave, and thus, the classic subgradient method to solve the dual problem is given by
\begin{equation}\label{eq:Lambda-Update}
    \boldsymbol{\lambda}^{(k+1)} = \boldsymbol{\lambda}^{(k)} + \gamma_k\mathbf{d}^{{(k)}},
\end{equation}
where $\gamma_k>0$ is the stepsize and $\mathbf{d}^{(k)}$ is a supergradient of $g$ at $\boldsymbol{\lambda}^{(k)}\in\R^{n}$, and $k$ signifies the iteration index. Note that the subproblem coordination is solely required to jointly construct the subgradient $\mathbf{d}^{(k)}$ at iterate $k$. The resulting algorithm has a very simple form.
\begin{breakablealgorithm}
	\caption{Dual Decomposition Algorithm} 
	\begin{algorithmic}[1]
	    \Require $\boldsymbol{\lambda}^{(0)}\in \R^{n(m-1)}$.
	    \State $k=0$.
		\Repeat 
            \State Solve subproblems  with $\boldsymbol{\lambda}{=}\boldsymbol{\lambda}^{(k)}$ to yield $\mathbf{y}^{(k)}\in\R^{nm}$. 
            \State Compute $\mathbf{d}^{{(k)}}{=}\big[\big({\mathbf{y}}^{{(k)}}_1-{\mathbf{y}}^{{(k)}}_2\big)\tran \ \ldots \  \big({\mathbf{y}}^{{(k)}}_{{m-1}}-{\mathbf{y}}^{{(k)}}_{m}\big)\tran\big]\tran$. 
            \State  $\boldsymbol{\lambda}^{(k+1)} = \boldsymbol{\lambda}^{(k)} + \gamma_k\mathbf{d}^{(k)}$.
            \State $k:=k+1$.
        \Until{a stopping criterion true}
	\end{algorithmic} 
	\label{Alg:Base-Line}
\end{breakablealgorithm}

Under mild technical conditions, the convergence of $\boldsymbol{\lambda}^{(k)}$ to the solution $\boldsymbol{\lambda}^\star$ of \eqref{eq:dual-problem} can be ensured, i.e., $\boldsymbol{\lambda}^{(k)}\rightarrow \boldsymbol{\lambda}^\star$ \cite[\S~1.4.2]{Polyak-Intro-Opt-1987}. Thus, together with additional assumptions such as the strong duality between \eqref{eq:distributed-problem} and \eqref{eq:dual-problem}, at the termination of the algorithm, a reasonable guess for the solution $\mathbf{y}^{\star}$ of \eqref{eq:distributed-problem} is obtained by averaging the solutions of the subproblems, i.e., $(1/m)\sum_i \mathbf{y}^{(k)}_i$ \cite[\S~5.5.5]{Boyd-Vandenberghe-2004}.

It is worth pointing out that, alternatively, one may rely on other classic approaches, such as the alternating direction method of multipliers~\cite{Boyd-Parikh-Chu-Peleato-Eckstein-2010} and proximal gradient methods\cite{Boyd-Parikh-Proximal-Algorithms-2014}, for developing distributed solution methods for problem~\eqref{eq:main-problem}. Analysis of such methods under inexact settings is of interest in its own~right and extraneous to the main focus of this paper.

\section{Inexact Algorithms}\label{sec:Imperfect-Coordination}


We consider the case where the subproblem coordination in each iteration $k$ is not perfect, \cf line 4 of Algorithm~\ref{Alg:Base-Line}. In particular, instead of the exact $\mathbf{y}_i^{(k)}$, a distorted vector $\hat {\mathbf{y}}_i^{(k)}$ is used when computing $\mathbf{d}^{(k)}$. As a result, instead of the exact $\mathbf{d}^{(k)}$, a distorted vector $\hat {\mathbf{d}}^{(k)}$ given by
\begin{equation}\label{eq:Distorted-Subgradient}
    \hat {\mathbf{d}}^{(k)}=\big[\big(\hat{\mathbf{y}}_1^{(k)}-\hat{\mathbf{y}}_2^{(k)})\tran \ \ldots \ (\hat{\mathbf{y}}_{m-1}^{(k)}-\hat{\mathbf{y}}_m^{(k)}\big)\tran\big]\tran
\end{equation}
is used in the dual variable update of Algorithm~\ref{Alg:Base-Line}, \cf line 5.

The distortion associated with $\mathbf{y}_i^{(k)}$, $i=1,\ldots,m$, $k\in\mathbb{Z}^0_+$ is denoted by $\mathbf{r}_i^{(k)}\in\R^n$, where
\begin{equation}\label{eq:Distortion}
    \hat {\mathbf{y}}_i^{(k)} = \mathbf{y}_i^{(k)} + \mathbf{r}_i^{(k)}.
\end{equation}
The distortion $\mathbf{r}_i^{(k)}$ can model numerous inexact settings, as remarked below.

\begin{remark}\label{Remark:Genarality-of-Error}
The \emph{additive} distortion $\mathbf{r}_i^{(k)}$ can model errors in many large-scale optimization problems, including quantization errors \cite{sindri_Communication_Complexity_2018,Sindri-2018,Sindri-2020}, approximation errors \cite{hu-biased-SGD-2020,Chen-biased-SGD-2018}, errors due to subproblem solver accuracy \cite{Yifan-inexact-dual-decomposition-2021,Necoara-inexact-dual-decomposition-2014}, errors in dual variable coordination, noise-induced in wireless settings \cite{Carlo-Machine-Learning}, and measurement errors, among others.
\end{remark}

It is worth pointing out that distortions in dual variables can be modeled as a distortion of underlying primal variables in general. To this end, one has to appeal to classic theories such as perturbation and sensitivity analysis~\cite{Anthony-Purturb-1983,Giorgio-Purturb-2018}. Our results hold as long as such dual-domain distortions are reflected in primal-domain, as in \eqref{eq:Distortion}. Despite the generality of $\mathbf{r}_i^{(k)}$, we refer to it as distortion due to imperfect coordination unless otherwise specified. When modeling the distortion $\mathbf{r}^{(k)}_i$, we assume nothing except the norm boundedness of the distortion. More specifically, we have the following assumption about the distortion $\mathbf{r}_i^{(k)}$.
\begin{assump}[Absolute Deterministic Distortion]\label{Assumption:Absolute-Deterministic-Error}
The distortion $\mathbf{r}_i^{(k)}$ associated with $\mathbf{y}_i^{(k)}$, $i=1,\ldots,m$, $k\in\mathbb{Z}^0_+$ is bounded by $\varepsilon_i\in\R$, i.e.,
\begin{equation}\label{eq:Absolute-Deterministic-Error}
  \|\mathbf{r}_i^{(k)}\|\leq \varepsilon_i, \quad i=1,\ldots,m, \quad k\in\mathbb{Z}^0_+.  
\end{equation}
\end{assump}
That is, the coordination of $\mathbf{y}_i^{(k)}$, $i=1,\ldots,m$, $k\in\mathbb{Z}^0_+$, always undergoes a persistent distortion \emph{bounded} by a worst-case value $\varepsilon_i$. Note that the worst-case characteristics of the distortion remain intact despite the iteration number $k$. Moreover, the distortion need not be random. If it is random, it need not be stationary, uncorrelated, or even zero mean. Roughly speaking, the considered distortions can be thought of as if they are imposed by some physical phenomena or by systems' imperfections that are unavoidable~\cite{Alessandri-bounded-error-2002}.

\begin{figure}[t!]
\centering
\begin{tikzpicture}[scale=0.5]
\tikzset{vertex/.style = {shape=circle,draw,minimum size=0.2 cm}}
\tikzset{edge/.style = {->,> = latex}}
\node[vertex] (1) at  (0,0) {CN};
\node[vertex] (2) at  (-7,-1.5) {$1$};
\node[vertex] (3) at  (-2,-1.5) {$2$};
\node[vertex] (4) at  (2,-1.5) {$3$};
\node[vertex] (5) at  (7,-1.5) {$m$};
 \draw [thick,stealth-]  (1.175) -- node[left, xshift=1.3cm, yshift=0.6cm] {$\hat{\mathbf{y}}_1^{(k)}$} (2.30);
 \draw [thick,line width=0.5mm,-stealth]  (1.190) -- node[right, xshift=-1.4cm, yshift=-0.6cm] {$\boldsymbol{\lambda}_1^{(k+1)}$} (2.5);
\draw [thick,stealth-]  (1.230) -- (3.60);
\draw [thick, dotted]  (2.9,-1.5) -- (5.7,-1.5);
\draw [thick,line width=0.5mm,-stealth]  (1.250) -- (3.35);
\draw [thick,stealth-]  (1.310) -- (4.120);
\draw [thick,line width=0.5mm,-stealth]  (1.290) -- (4.145);
\draw [thick,stealth-]  (1.5) -- node[right, xshift=-1.2cm, yshift=0.6cm] {$\hat{\mathbf{y}}_m^{(k)}$} (5.150);
\draw [thick,line width=0.5mm,-stealth]  (1.350) -- node[left, xshift=1.4cm, yshift=-0.7cm] {$\boldsymbol{\lambda}_{m-1}^{(k+1)}$} (5.170);
\end{tikzpicture}
\caption{Communication Structure: Partially Distributed Algorithm.} 
\label{Fig:Partially-Distributed-Algorithm}
\end{figure}
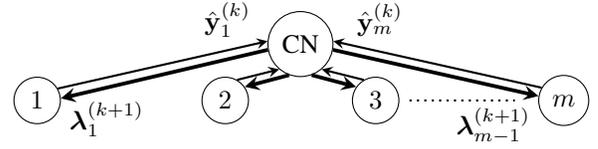

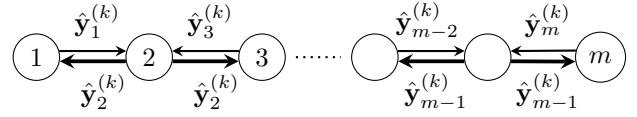
\begin{figure}[t!]
\centering
\begin{tikzpicture}[scale=0.5]
\tikzset{vertex/.style = {shape=circle,draw,minimum size=0.2cm}}
\tikzset{edge/.style = {->,> = latex}}
\node[vertex] (1) at  (0,0) {$1$};
\node[vertex] (2) at  (3,0) {$2$};
\node[vertex] (3) at  (6,0) {$3$};
\node[vertex] (4) at  (9,0) {\addw{$1$}};
\node[vertex] (5) at  (12,0) {\addw{$1$}};
\node[vertex] (6) at  (15,0) {$m$};
\draw [thick, dotted]  (6.9,0) -- (8.1,0);
 \draw [thick,-stealth]  (1.15) -- node[left, xshift=0.5cm, yshift=0.4cm] {$\hat{\mathbf{y}}_1^{(k)}$} (2.165);
 \draw [thick,line width=0.5mm,stealth-]  (1.350) -- node[right, xshift=-0.3cm, yshift=-0.4cm] {$\hat{\mathbf{y}}_2^{(k)}$} (2.190);
 \draw [thick, stealth-]  (2.15) -- node[left, xshift=0.5cm, yshift=0.4cm] {$\hat{\mathbf{y}}_3^{(k)}$} (3.165);
 \draw [thick,line width=0.5mm,-stealth]  (2.350) -- node[right, xshift=-0.3cm, yshift=-0.4cm] {$\hat{\mathbf{y}}_2^{(k)}$} (3.190);
 \draw [thick,-stealth]  (4.15) -- node[left, xshift=0.5cm, yshift=0.4cm] {$\hat{\mathbf{y}}_{m-2}^{(k)}$} (5.165);
 \draw [thick,line width=0.5mm,stealth-]  (4.350) -- node[right, xshift=-0.5cm, yshift=-0.4cm] {$\hat{\mathbf{y}}_{m-1}^{(k)}$} (5.190);
 \draw [thick,stealth-]  (5.15) -- node[left, xshift=0.5cm, yshift=0.4cm] {$\hat{\mathbf{y}}_{m}^{(k)}$} (6.165);
 \draw [thick,line width=0.5mm,-stealth]  (5.350) -- node[right, xshift=-0.5cm, yshift=-0.4cm] {$\hat{\mathbf{y}}_{m-1}^{(k)}$} (6.190);
\end{tikzpicture}
\caption{Communication Structure: Fully Distributed Algorithm.}
		\label{Fig:Fully-Distributed-Algorithm}
	\end{figure}

Our exposition of imperfect coordination can be centered on different variants of Algorithm~\ref{Alg:Base-Line}. A setting where a central node (CN) solely performs the subproblem coordination is considered as a partially distributed variant, \cf~\figurename~\ref{Fig:Partially-Distributed-Algorithm}. Federated learning~\cite{Federated-learning-McMahan-2016} which has received significant attention in recent literature, closely falls under this setting. 
In addition, a fully distributed variant, in the sense that there is no central authority, can also be considered, \cf~\figurename~\ref{Fig:Fully-Distributed-Algorithm}. Without loss of generality from a technical standpoint, let us consider the latter in our subsequent exposition.

It turns out that the decomposition structure [\cf \figurename~\ref{Fig:Decomposition-Structure}] considered when reformulating problem~\eqref{eq:distributed-problem} suggests a subproblem coordination mechanism where only the communication between neighboring SSs is necessary. The communication structure is depicted in \figurename~\ref{Fig:Fully-Distributed-Algorithm}. Associated with each $i\in \{1,\ldots,m\}$, the following resources will enable the subproblem coordination: 1) an error-free communication channel from SS $i{-}1$ to $i$, 2) an error-free communication channel from SS $i{+}1$ to $i$.

In the case above, the distortion $\mathbf{r}_i^{(k)}$ is due to inevitable approximation errors that come about
as a result of attempts to reduce communication overhead, see \cite{Sindri-2020,sindri_Communication_Complexity_2018}.
The resulting algorithm is summarized below.

\begin{breakablealgorithm}
	\caption{Fully Distributed Algorithm} 
	\begin{algorithmic}[1]
	    \Require $\boldsymbol{\lambda}^{(0)}\in\R^{n(m-1)}$; $\boldsymbol{\lambda}_0^{(j)}=\boldsymbol{\lambda}_m^{(j)}=\mathbf{0}\in\R^n, \ j\in\mathbb{Z}_+^0$; $\mathbf{y}_0^{(j)}=\mathbf{y}_{m+1}^{(j)}=\mathbf{0}\in\R^n, \ j\in\mathbb{Z}_+^0$.
	    \State $k=0$.
		\Repeat 
            \State $\forall i$, SS $i$ computes $\mathbf{y}_i^{(k)}$ by solving 
\begin{equation} \nonumber
\begin{array}{ll}
{\mbox{minimize}}_{\mathbf{y}_i\in\mathcal{Y}} & f_i(\mathbf{y}_i)+ {\big(\boldsymbol{\lambda}^{(k)}_i-\boldsymbol{\lambda}^{(k)}_{i-1}\big)}\tran \mathbf{y}_i. 
\end{array}
\end{equation}
\State  $\forall i$, SS  $i$ transmits $\hat{\mathbf{y}}_i^{(k)}$ to $i-1$ and $i+1$, \cf \eqref{eq:Distortion}. 
            \State  $\forall i$, SS $i$ receives $\hat{\mathbf{y}}_{i-1}^{(k)}$ and $\hat{\mathbf{y}}_{i+1}^{(k)}$ from $i{-}1$ and $i{+}1$.
\State $\forall i$, SS $i$ sets $\hat {\mathbf{d}}_i^{(k)}{=}\big[(\hat{\mathbf{y}}_{i-1}^{(k)}-\hat{\mathbf{y}}_{i}^{(k)})\tran \ (\hat{\mathbf{y}}_{i}^{(k)}{-}\hat{\mathbf{y}}_{i{+}1}^{(k)})\tran\big]\tran$. 
            \State $\forall i$, SS $i$ sets
            \begin{equation*}
             \hspace{2.5mm}\big[(\boldsymbol{\lambda}^{(k{+}1)}_{i{-}1})\tran \ (\boldsymbol{\lambda}^{(k{+}1)}_{i})\tran\big]\tran{=}\big[(\boldsymbol{\lambda}^{(k)}_{i{-}1})\tran \ (\boldsymbol{\lambda}^{(k)}_{i})\tran\big]\tran{+}\gamma_k\hat {\mathbf{d}}^{(k)}_i
            \end{equation*}
            \State $k:=k+1$
        \Until{a stopping criterion true}
	\end{algorithmic} 
	\label{Alg:Fully-Distributed}
\end{breakablealgorithm}

Each operation listed from $3{-}7$ is conducted in parallel among $i\in\{1,\ldots, m\}$. Here the subproblem coordination is solely achieved by each SS's communication with its neighbors and local computations, \cf lines 4-6 of Algorithm~\ref{Alg:Fully-Distributed}. Subproblem coordination stage 1 [see line 4] is the source of the imperfect coordination, where an absolute deterministic distortion is introduced [\cf Assumption~\ref{Assumption:Absolute-Deterministic-Error}].

Attempts at ascertaining convergence of iterative algorithms always rely on certain key characteristics of the underlying functions. In this respect, related to Algorithm~\ref{Alg:Fully-Distributed}, characteristics of the primal function $f$ play a crucial role. We now proceed toward hypothesizing some such characteristics. 

\begin{assump}[Strong Duality]\label{Assumption:Strong-Duality}
The optimal values $p^\star$ and $d^\star$ of problems \eqref{eq:distributed-problem} and \eqref{eq:dual-problem}, respectively, are attained. Moreover, the strong duality between \eqref{eq:distributed-problem} and \eqref{eq:dual-problem} holds, i.e.,
\begin{equation}
     p^\star= f(\mathbf{y}^\star)= g(\boldsymbol{\lambda}^\star)=d^\star,
\end{equation}
for some $\mathbf{y}^\star\in\{\mathbf{y}\in\R^{nm} \ | \ \forall~i \  \mathbf{y}_i\in\mathcal{Y}, \  \mathbf{A}\mathbf{y}=\mathbf{0} \}$ and for some $\boldsymbol{\lambda}^\star\in\R^{n(m-1)}$, where $\mathbf{A}$ is defined in~\eqref{eq:A-Matrix-for-Consensus-Constraint}.
\end{assump}

Assumption~\ref{Assumption:Strong-Duality} bridges the primal and dual domains through which the dual convergence results are maneuvered in deriving primal convergence results, for otherwise not viable in general. The assumption usually holds in general for convex problems under some constraint qualifications~\cite[\S~5.2.3]{Boyd-Vandenberghe-2004}.
\begin{assump2}[Strongly Convex Objectives at SSs]\label{Assumption:Strongly-Convex-f}
 The functions $f_i$s in problem~\eqref{eq:distributed-problem} are strongly convex with constant $\mu_i>0$, $i=1,\ldots,m$.
\end{assump2}
\begin{assump2}[Gradient Lipschitz Continuous Objectives at SSs]\label{Assumption:Lipschitz-Convex-Grad-f}
 The set $\mathcal{Y}$ in problem~\eqref{eq:distributed-problem} equals $\R^n$. Moreover, $f_i$s are differentiable, and the gradients $\nabla f_i$s are Lipschitz continuous on $\R^n$ with constant $L_i>0$, $i=1,\ldots,m$. 
\end{assump2}
Note that Assumption~\ref{Assumption:Strongly-Convex-f} usually holds in many application domains because the primal objective functions can inherently be strongly convex~\cite{Rabbat-strong-convexity-2012}. Assumption~\ref{Assumption:Lipschitz-Convex-Grad-f} can be interpreted as a bound on the
second derivative of primal objective functions, which is usually true in many practical problems, e.g., \cite{Jerome-logistic-2010}. Important consequences of those assumptions, which have been heavily used in the sequel, are deferred to Appendix~\ref{sec:Characteristics-of-Dual-Function}.

\section{Convergence Analysis}\label{sec:Convergence-Analysis-I}

Outfitted with the Lipschitzian and strong convexity properties of $-g$, let us now derive the convergence properties of Algorithm~	\ref{Alg:Fully-Distributed}. 
Roughly speaking, the exposition in this section is divided into two main cases: 
\begin{align}
    & \mbox{\textsc{Case 1}: Assumption~\ref{Assumption:Strongly-Convex-f} holds, \cf Proposition~\ref{Proposition:Lipscitz-Continuity-of-Grad-g}.} \label{eq:CASE-1}\\
    & \mbox{\textsc{Case 2}: Assumption~\ref{Assumption:Strongly-Convex-f} and Assumption~\ref{Assumption:Lipschitz-Convex-Grad-f} hold},\nonumber \\
    &\quad\quad \quad \hspace{3mm} \mbox{\cf Proposition~\ref{Proposition:Lipscitz-Continuity-of-Grad-g} and Proposition~\ref{Proposition:Strong-Concavity-of-g}.} \label{eq:CASE-2}
    \end{align}
For each case, it is useful to restate the Lagrange multiplier update performed by Algorithm~\ref{Alg:Fully-Distributed}, i.e.,
\begin{equation}\label{eq:Lambda-Update-Distorted}
    \boldsymbol{\lambda}^{(k+1)} = \boldsymbol{\lambda}^{(k)} + \gamma_k\hat {\mathbf{d}}^{(k)},
\end{equation}
see line 7 of Algorithm~\ref{Alg:Fully-Distributed}. Moreover, recall that the primal solution computed by Algorithm~\ref{Alg:Fully-Distributed} in each iteration $k\in\mathbb{Z}_+^0$ is $\mathbf{y}^{(k)}$, where 
\begin{equation}\label{eq:Local-y-k-in-vector-form}
   \hspace{-1mm} \mathbf{y}^{(k)}{=}\big[(\mathbf{y}^{(k)}_1)\tran \ {\ldots} \ (\mathbf{y}^{(k)}_m)\tran\big]\tran{=}\underset{\mathbf{y}\in\bar{\mathcal{Y}}}{\arg\min} \  f(\mathbf{y}){+} (\boldsymbol{\lambda}^{(k)})\tran \mathbf{A} \mathbf{y}.
\end{equation}

We have derived convergence results of the sequences of Lagrange multipliers $\{\boldsymbol{\lambda}^{(k)}\}$, the primal solutions $\{\mathbf{y}^{(k)}\}$, and the primal objective function values $\{f(\mathbf{y}^{(k)})\}$. Then we discuss how a feasible point $\tilde {\mathbf{y}}^{(k)}$ to problem~\eqref{eq:distributed-problem} is computed by using $\mathbf{y}^{(k)}$. More importantly, convergences of the sequences $\{\tilde {\mathbf{y}}^{(k)}\}$ and $\{f(\tilde {\mathbf{y}}^{(k)})\}$ are also mathematically substantiated. 

All derivations are conducted under two stepsize rules: 1) constant stepsize rule, i.e., $\gamma_k=\gamma, \ \forall k$; 2) nonsummable stepsize rule, i.e., $\sum_{k=0}^{\infty}\gamma_k=\infty$.
Note that the constant stepsize rule is a particular case of the nonsummable stepsize rule.



\subsection{Key Remarks and Related Results}
Let us start by outlining some results that are useful in the latter part of the section. A consequence of Assumption~\ref{Assumption:Absolute-Deterministic-Error} is an upperbound on the overall distortion due to imperfect subproblem coordination, which is outlined below.

\begin{remark}\label{Remark:Assumption:Absolute-Deterministic-Error-All}
Assumption~\ref{Assumption:Absolute-Deterministic-Error} entails an absolute deterministic distortion of $\mathbf{d}^{(k)}$ [\cf \eqref{eq:Lambda-Update} and \eqref{eq:Distorted-Subgradient}]. In particular, we have,
\begin{equation} \label{eq:total-error-bound}
  \|\hat {\mathbf{d}}^{(k)}-\mathbf{d}^{(k)}\|\leq \epsilon
\end{equation}
for Algorithm~\ref{Alg:Fully-Distributed}, where  $\epsilon=\sqrt{\sum_{i=1}^{m-1}(\varepsilon_i+\varepsilon_{i+1})^2}$. 
\end{remark}

Let $h=-g$ for clarity. It is worth noting that the function $h$ is differentiable, as remarked below.
\begin{remark}\label{Remark:Differentiability-of-h}
Assumption~\ref{Assumption:Strongly-Convex-f} entails the differentiability of $h$ on $\R^{n(m-1)}$.
\end{remark}

Finally, we record a lemma highlighting a recursive inequality that is useful when asserting convergence of both cases, i.e., \textsc{Case 1} and \textsc{Case 2}.

\begin{lemma} \label{Lemma:II-General-Recursive-Inequality-For-Both-Cases}
Suppose Assumption~\ref{Assumption:Absolute-Deterministic-Error} 
and Assumption~\ref{Assumption:Strongly-Convex-f} hold. Let $\gamma_k$ satisfy the condition $0<\gamma_k\leq 1/L_h$ for all $k\in\mathbb{Z}^0_+$. Then, the function $h$ evaluated at Lagrange multipliers computed in consecutive iterations $k$, and $k+1$ of Algorithm~\ref{Alg:Fully-Distributed} conforms~to
 %
 \begin{equation}
 h(\boldsymbol{\lambda}^{(k+1)})\leq  h(\boldsymbol{\lambda}^{(k)})-\frac{\gamma_k}{2}\|\nabla h(\boldsymbol{\lambda}^{(k)})\|^2+\frac{\gamma_k}{2}\epsilon^2,
   \label{eq:Recursive-Inequality-For-Both-Cases}
 \end{equation}
 where $L_h=(1/\mu)\left(2+2\cos(\pi/m)\right)$, with $\mu=\min_i~\mu_i$. 
\end{lemma}

\begin{IEEEproof}
Assumption~\ref{Assumption:Strongly-Convex-f} entails that
Proposition~\ref{Proposition:Lipscitz-Continuity-of-Grad-g} holds. Use descent lemma \cite[\S~5.1.1]{Amir-first-order-methods-2017} to yield 
\begin{equation} \nonumber
h(\boldsymbol{\gamma})\leq h(\boldsymbol{\delta})+ \nabla h(\boldsymbol{\delta})\tran(\boldsymbol{\gamma}-\boldsymbol{\delta})+\frac{L_h}{2}\|\boldsymbol{\gamma}-\boldsymbol{\delta}\|^2~
\end{equation}
for all  $\boldsymbol{\gamma},\boldsymbol{\delta}\in\R^{n(m-1)}$. We start by setting $\boldsymbol{\gamma}=\boldsymbol{\lambda}^{(k+1)}$ and $\boldsymbol{\delta}=\boldsymbol{\lambda}^{(k)}$. Thus, the result follows by noting $\boldsymbol{\gamma}-\boldsymbol{\delta}=\gamma_k\hat {\mathbf{d}}^{(k)}$ [\cf \eqref{eq:Lambda-Update-Distorted}], $\mathbf{r}^{(k)}=\hat {\mathbf{d}}^{(k)}-\mathbf{d}^{(k)}$, $\mathbf{d}^{(k)}=-\nabla h(\boldsymbol{\lambda}^{(k)})$ [\cf Remark~\ref{Remark:Differentiability-of-h}], together by restricting $0<\gamma_k\leq 1/L_h$ for all $k\in\mathbb{Z}^0_+$ and using that $\|\mathbf{r}^{(k)}\|\leq \epsilon$ [\cf~Remark~\ref{Remark:Assumption:Absolute-Deterministic-Error-All}].
\end{IEEEproof}
%


\subsection{Case 1: Algorithm~\ref{Alg:Fully-Distributed} under Assumption~\ref{Assumption:Strongly-Convex-f}}\label{subsec:Lipschitz-continuous}

In this section, we derive convergence results in the dual-domain and those related to primal optimality violations for Algorithm~\ref{Alg:Fully-Distributed}, along with their rates of convergence. 



\begin{lemma} \label{Lemma:Lipschitz-continuous-constant-nonsum}
Suppose Assumption~\ref{Assumption:Absolute-Deterministic-Error} 
and Assumption~\ref{Assumption:Strongly-Convex-f} hold. Let $\{\boldsymbol{\lambda}^{(k)}\}$ be the sequence of Lagrange multipliers generated by Algorithm~\ref{Alg:Fully-Distributed} with stepsize $\gamma_k$ conforming to $0<\gamma_k\leq1/L_h$ for all $k\in\mathbb{Z}^0_+$. Then 
\begin{equation}  \label{eq:bound-min-normgrad-square-nonsum}
\underset{i\in\{0,\ldots,k\}}{\min} \ \|\nabla h(\boldsymbol{\lambda}^{(i)})\| \leq \displaystyle \sqrt{\frac{2\left(h(\boldsymbol{\lambda}^{(0)})-h(\boldsymbol{\lambda}^\star)\right)}{\sum_{i=0}^{k}\gamma_i}}+\epsilon,
\end{equation}
where $L_h=(1/\mu)\left(2+2\cos(\pi/m)\right)$, with $\mu=\min_i~\mu_i$.
\end{lemma}
\begin{IEEEproof}
By applying recursively \eqref{eq:Recursive-Inequality-For-Both-Cases} [\cf Lemma~\ref{Lemma:II-General-Recursive-Inequality-For-Both-Cases}] and rearranging the resulting terms yields
\begin{align} \nonumber
\sum_{i=0}^{k}\gamma_i \|\nabla h(\boldsymbol{\lambda}^{(i)})\|^2&\leq 2(h(\boldsymbol{\lambda}^{(0)})-h(\boldsymbol{\lambda}^{(k+1)}))+\epsilon^2  \sum_{i=0}^k \gamma_i \allowdisplaybreaks\\ \label{eq:Case-I-Recursion-constant-nonsum-2-in}
&\leq 2\left(h(\boldsymbol{\lambda}^{(0)})-h(\boldsymbol{\lambda}^\star)\right)+\epsilon^2  \sum_{i=0}^k \gamma_i, \allowdisplaybreaks
\end{align}
where the last inequality is immediate from that $h(\boldsymbol{\lambda}^{(k)})\geq h(\boldsymbol{\lambda}^{\star})$, for all $k\in\mathbb{Z}^0_+$ and $\boldsymbol{\lambda}^\star$ is a dual solution [\cf Assumption~\ref{Assumption:Strong-Duality}]. Then the result follows by further rearranging the terms of \eqref{eq:Case-I-Recursion-constant-nonsum-2-in}, together with that ${\min}_j \ \|\nabla h(\boldsymbol{\lambda}^{(j)})\|^2\leq \|\nabla h(\boldsymbol{\lambda}^{(i)})\|^2$ for all $i\in\{0,\ldots,k\}$, $\sqrt{\min_j\|\cdot\|^2}=\min_j\sqrt{\|\cdot\|^2}$, and $\sqrt{x+y}\leq\sqrt{x}+\sqrt{y}$ for all $x,y\geq 0$. 
\end{IEEEproof}

\begin{corollary} \label{Corollary:Lipschitz-continuous-nonsum-step}
Suppose Assumption~\ref{Assumption:Absolute-Deterministic-Error}
and Assumption~\ref{Assumption:Strongly-Convex-f} hold. Let $\{\boldsymbol{\lambda}^{(k)}\}$ be the sequence of Lagrange multipliers generated by Algorithm~\ref{Alg:Fully-Distributed}. Moreover, suppose $\gamma_k$ satisfy the nonsummable  stepsize rule with $0<\gamma_k\leq1/L_h$. Then
\begin{enumerate}
    \item $\displaystyle \underset{k}{{\limsup}} \ \underset{i\in\{0,\ldots,k\}}{\min} \ \|\nabla h(\boldsymbol{\lambda}^{(i)})\| \leq \epsilon$.
    \item for $\gamma_k=\gamma/(k+1)^p$, where $0<\gamma\leq1/L_h$ and $0\leq p\leq1$,
%
\begin{equation*}\label{eq:Lipschitz-continuous-nonsum-Step-convergence-rate}
\underset{i\in\{0,\ldots,k\}}{\min} \ \|\nabla h(\boldsymbol{\lambda}^{(i)})\| {=}
    \begin{cases}
    \displaystyle O\left(\frac{1}{\sqrt{k^{{1-p}}}}\right)+\epsilon, & \text{$p\in[0,1)$}\\
         \displaystyle O\left( \frac{1}{\sqrt{\log k}} \right) +\epsilon, & \text{$p=1$},
    \end{cases}       
\end{equation*}
and the best convergence rate is of the order $O(1/\sqrt{k})$, which is achieved when $p=0$, i.e., the fixed stepsize.
\end{enumerate}
\end{corollary}
\begin{IEEEproof}
This is clear from Lemma~\ref{Lemma:Lipschitz-continuous-constant-nonsum} and straightforward algebraic manipulations and is omitted.
\end{IEEEproof}

Corollary~\ref{Corollary:Lipschitz-continuous-nonsum-step} indicates that the least upperbound of $\min_i\|\nabla(h(\boldsymbol{\lambda}^{(i)})\|$ can converges to a neighborhood around~$0$ at a rate of order $O(1/\sqrt{k})$, where the size of the neighborhood depends on $\epsilon$, \cf \eqref{eq:total-error-bound}. It is straightforward to see that the fastest rate corresponds to the fixed stepsize rule with $\gamma_k=1/L_h$ for all $k\in\mathbb{Z}^0_+$, \cf \eqref{eq:bound-min-normgrad-square-nonsum}. We note that under Assumption~\ref{Assumption:Strongly-Convex-f} and ideal conditions, the standard gradient descent on the dual function can converge to the optimal value at a rate of order $O(1/k)$~\cite[Theorem~2.1.14]{Nesterov-2018}. In contrast, under Assumption~\ref{Assumption:Strongly-Convex-f}, together with absolute deterministic errors, only a rate of ${O}(1/\sqrt{k})$ is established but with respect to the norm of the gradient. Our approach, however, restrict us to derive convergences with respect to function values. 

Convergence of the primal optimality violations is established next.

\begin{prop2}\label{Proposition:Lipschitz-continuous-g-Primal-Result}
Suppose Assumption~\ref{Assumption:Absolute-Deterministic-Error}, Assumption~\ref{Assumption:Strong-Duality}, and Assumption~\ref{Assumption:Strongly-Convex-f} hold. Let $\{\boldsymbol{\lambda}^{(k)}\}$ be the sequence of Lagrange multipliers generated by Algorithm~\ref{Alg:Fully-Distributed} and $\{\mathbf{y}^{(k)}\}$ be the corresponding sequence of primal variables. Moreover, suppose that the functions $f_i$, $i=1,\ldots,m$ are differentiable. Let $\gamma_k$ satisfy the nonsummable  stepsize rule with $0<\gamma_k\leq1/L_h$. If the distance to the dual optimal solution $\|\boldsymbol{\lambda}^{(k)}-\boldsymbol{\lambda}^\star\|$ is uniformly bounded by some scalar $D$, then
\begin{enumerate}
    \item $\displaystyle \underset{k}{{\limsup}} \ \underset{i\in\{0,\ldots,k\}}{\min} \ \|\mathbf{y}^{(i)}-\mathbf{y}^\star\| \leq \sqrt{2D\epsilon/\mu}$.
\item  $\displaystyle \underset{k}{{\limsup}} \ \underset{i\in\{0,\ldots,k\}}{\min} \ \big(f(\mathbf{y}^{(i)})-f(\mathbf{y}^\star)\big) \leq D\epsilon+\sqrt{D}S(D+\|\boldsymbol{\lambda}^{\star}\|) \sqrt{\epsilon}$, where the positive scalar $S=$ $\sqrt{(4+4\cos(\pi/m))/\mu}$.
    \item for $\gamma_k=\gamma/(k+1)^p$, where $0<\gamma\leq1/L_h$ and $0\leq p\leq1$,
\begin{equation*}\label{eq:Lipschitz-continuous-nonsum-Step-convergence-rate}
\underset{i\in\{0,\ldots,k\}}{\min} \ \|\mathbf{y}^{(i)}-\mathbf{y}^\star\| =
    \begin{cases}
    \displaystyle O\left(\frac{1}{\sqrt[4]{k^{{1-p}}}}\right)+\sqrt{{2D\epsilon}/{\mu}}, &\\
    &\hspace{-12mm}\text{$p\in[0,1)$}\\
         \displaystyle O\left( \frac{1}{\sqrt[4]{\log k}} \right) +\sqrt{{2D\epsilon}/{\mu}}, \\
         &\hspace{-8mm} \text{$p=1$},
    \end{cases}       
\end{equation*}
and the best convergence rate is of the order $O(1/\sqrt[4]{k})$, which is achieved when $p=0$.
 \item for $\gamma_k=\gamma/(k+1)^p$, where $0<\gamma\leq1/L_h$ and $0\leq p\leq1$, 
\begin{equation*}\label{eq:Lipschitz-continuous-nonsum-Step-convergence-rate}
\underset{i\in\{0,\ldots,k\}}{\min} \ \|f(\mathbf{y}^{(i)})-f(\mathbf{y}^\star)\| =
    \begin{cases}
    \displaystyle O\left(\frac{1}{\sqrt[4]{k^{{1-p}}}}\right)+D\epsilon \\ +\sqrt{D}S(D+\|\boldsymbol{\lambda}^{\star}\|) \sqrt{\epsilon},& \\
    & \hspace{-17mm}\text{$p\in[0,1)$}\\
         \displaystyle O\left( \frac{1}{\sqrt[4]{\log~k}} \right) +D\epsilon\\
         +\sqrt{D}S(D+\|\boldsymbol{\lambda}^{\star}\|) \sqrt{\epsilon}, & \\
         & \hspace{-13mm}\text{$p=1$},
    \end{cases}       
\end{equation*}
and the best convergence rate is of the order $O(1/\sqrt[4]{k})$, which is achieved when $p=0$.
\end{enumerate}
\end{prop2}
\begin{IEEEproof}
See Appendix~\ref{Appendix:Proposition:Lipschitz-continuous-g-Primal-Result}.
\end{IEEEproof}
Note that Proposition~\ref{Proposition:Lipschitz-continuous-g-Primal-Result} for the primal-domain convergences relies on a few additional hypotheses, unlike Corollary~\ref{Corollary:Lipschitz-continuous-nonsum-step}, which could not be avoided in our approach to proving the results.

\subsection{Case 2: Algorithm~\ref{Alg:Fully-Distributed} under Assumption~\ref{Assumption:Strongly-Convex-f} and Assumption~\ref{Assumption:Lipschitz-Convex-Grad-f}}\label{subsec:Strongly-convex-Lipschits-continuous}

\begin{lemma} \label{Lemma:Strongly-convex-Lipschitz-continuous-constant-nonsum-step}
Suppose Assumption~\ref{Assumption:Absolute-Deterministic-Error},
Assumption~\ref{Assumption:Strongly-Convex-f}, and Assumption~\ref{Assumption:Lipschitz-Convex-Grad-f} hold. Let $\{\boldsymbol{\lambda}^{(k)}\}$ be the sequence of Lagrange multipliers generated by Algorithm~\ref{Alg:Fully-Distributed} with the stepsize $\gamma_k$ satisfying the condition $0<\gamma_k\leq1/L_h$ for all $k\in\mathbb{Z}^0_+$. Then
\begin{equation}  \label{eq:Strongly-convex-Lipschitz-continuous-constant-nonsum-step}
h(\boldsymbol{\lambda}^{(k+1)}){-}h(\boldsymbol{\lambda}^{\star}){\leq} (1{-}\gamma_k \mu_h)\left(h(\boldsymbol{\lambda}^{(k)}){-}h(\boldsymbol{\lambda}^{\star})\right){+}\frac{\gamma_k\epsilon^2}{2},
\end{equation}
where $L_h=(1/\mu)\left(2+2\cos(\pi/m)\right)$, with $\mu=\min_i~\mu_i$, and $\mu_h=(1/L)\left(2-2\cos(\pi/m)\right)$, with $L=\max_i~L_i$.
\end{lemma}
\begin{IEEEproof}
Since Assumption~\ref{Assumption:Lipschitz-Convex-Grad-f} holds, $h(\boldsymbol{\lambda})$ is strongly convex with constant $\mu_h$, \cf Proposition \ref{Proposition:Strong-Concavity-of-g}. Hence we have that $\|\nabla h(\boldsymbol{\lambda}^{(k)})\|^2\geq 2\mu_h\left(h(\boldsymbol{\lambda}^{(k)})-h(\boldsymbol{\lambda}^{\star})\right)$ \cite[p.~24]{Polyak-Intro-Opt-1987}. This, together with Lemma~\ref{Lemma:II-General-Recursive-Inequality-For-Both-Cases} yields
\begin{align}
    h(\boldsymbol{\lambda}^{(k+1)})-h(\boldsymbol{\lambda}^{\star}) \leq (1-\gamma_k\mu_h) \left(h(\boldsymbol{\lambda}^{(k)})-h(\boldsymbol{\lambda}^{\star})\right)+\frac{\gamma_k}{2}\epsilon^2\nonumber
\end{align}
the intended result.
\end{IEEEproof}
Note that we have $0\leq1-\gamma_k \mu_h<1$ in Lemma~\ref{Lemma:Strongly-convex-Lipschitz-continuous-constant-nonsum-step}, since $0<\gamma_k\leq 1/L_h$ and $\mu_h\leq L_h$. Let us next derive the convergence results for constant stepsize and for nonsummable stepsize. In the case of constant stepsize, the following result is immediate from Lemma~\ref{Lemma:Strongly-convex-Lipschitz-continuous-constant-nonsum-step}.

\begin{corollary} \label{Corollary:Strongly-convex-Lipschitz-continuous-constant-step}
Suppose Assumption~\ref{Assumption:Absolute-Deterministic-Error},
Assumption~\ref{Assumption:Strongly-Convex-f}, and Assumption~\ref{Assumption:Lipschitz-Convex-Grad-f} hold. Let $\{\boldsymbol{\lambda}^{(k)}\}$ be the sequence of Lagrange multipliers generated by Algorithm~\ref{Alg:Fully-Distributed} with the stepsize $\gamma_k=\gamma$ for all $k\in\mathbb{Z}^0_+$. 
Then for $0<\gamma\leq1/L_h$ 
\begin{equation}\label{eq:Strongly-convexL-ipschitz-continuous-constant-step}
\hspace{-2mm}    h(\boldsymbol{\lambda}^{(k)}){-}h(\boldsymbol{\lambda}^{\star})\leq  (1-\gamma\mu_h)^k \left(h(\boldsymbol{\lambda}^{(0)}){-}h(\boldsymbol{\lambda}^{\star})\right){+}\frac{\epsilon^2}{2\mu_h}. 
\end{equation}
Moreover, ${\limsup}_{k} \ \left(h(\boldsymbol{\lambda}^{(k)})-h(\boldsymbol{\lambda}^{\star})\right)\leq {\epsilon^2}/{(2\mu_h)}$.
\end{corollary}
\begin{IEEEproof}
The first part is immediate from the recursive application of \eqref{eq:Strongly-convex-Lipschitz-continuous-constant-nonsum-step} with $\gamma_k=\gamma$ and that $\sum_{i=0}^{k}(1-\gamma\mu_h)^i\leq 1/(\gamma\mu_h)$ for all $k\in\mathbb{Z}_+^0$. The second part follows because $\lim\sup_k (1-\gamma\mu_h)^k = 0$ and $h(\boldsymbol{\lambda}^{(0)})-h(\boldsymbol{\lambda}^{\star})<\infty$. 
\end{IEEEproof}

We note that \cite[Theorem~6]{Ahmad-biased-SGD-2021} is closely related to our result in a stochastic setting. According to Corollary~\ref{Corollary:Strongly-convex-Lipschitz-continuous-constant-step}, with constant stepsize, $h(\boldsymbol{\lambda}^{(k)})$ converges into a neighborhood of the optimal value $h(\boldsymbol{\lambda}^\star)$ with the rate of geometric progression, where the size of the neighborhood depends on $\epsilon$ [\cf \eqref{eq:total-error-bound}] and the constant $\mu_h$ that characterizes the strong convexity of $h$ [\cf Proposition~\ref{Proposition:Strong-Concavity-of-g}]. The following proposition asserts the convergences of the primal optimality violation.

\begin{prop2}\label{Proposition:Strongly-convex-Lipschitz-continuous-constant-step-Primal}
Suppose Assumption~\ref{Assumption:Absolute-Deterministic-Error}, Assumption~\ref{Assumption:Strong-Duality}, Assumption~\ref{Assumption:Strongly-Convex-f}, and Assumption~\ref{Assumption:Lipschitz-Convex-Grad-f} hold. Let $\{\boldsymbol{\lambda}^{(k)}\}$ be the sequence of Lagrange multipliers generated by Algorithm~\ref{Alg:Fully-Distributed} and $\{\mathbf{y}^{(k)}\}$ be the corresponding sequence of primal variables. Let the stepsize $\gamma_k=\gamma$ for all $k\in\mathbb{Z}^0_+$. 
Then for $0<\gamma\leq1/L_h$
\begin{enumerate}
    \item $\displaystyle \underset{k}{{\limsup}} \ \|\mathbf{y}^{(k)}-\mathbf{y}^\star\| \leq \epsilon\sqrt{1/(\mu\mu_h)}$.
\item  $\displaystyle \underset{k}{{\limsup}} \ \|f(\mathbf{y}^{(k)})-f(\mathbf{y}^\star)\| \leq  ({\epsilon^2}/{2\mu_h})\big(1+S\sqrt{2/\mu_h}\big) + \epsilon S\|\boldsymbol{\lambda}^{\star}\| \sqrt{1/{2\mu_h}}$, where the positive scalar $S=$ $\sqrt{(4+4\cos(\pi/m))/\mu}$.
    \item  the least upperbound of $ \|\mathbf{y}^{(k)}-\mathbf{y}^\star\|$ converges into a neighborhood of $0$ with the rate of geometric progression. 
 \item the least upperbound of $ \|f(\mathbf{y}^{(k)})-f(\mathbf{y}^\star)\|$ converges into a neighborhood of $0$ with the rate of geometric progression.
\end{enumerate}
\end{prop2}
\begin{IEEEproof}
The proof is virtually the same as that of  Proposition~\ref{Proposition:Lipschitz-continuous-g-Primal-Result}. Let $\mathbf{y}^{(k)}$ be the solution to problem~\eqref{eq:Compact-SubProblem-All-in-One}. Then
\begin{align}\label{eq:Lipschitz-continuous-Srong-Convex-g-Primal-Result-1-}
\hspace{-2mm}   \|\mathbf{y}^{(k)}{-}\mathbf{y}^\star\|^2 \leq  \frac{2}{\mu} (1{-}\gamma_k\mu_h)^k\left(h(\boldsymbol{\lambda}^{(0)}){-}h(\boldsymbol{\lambda}^{\star})\right){+}\frac{\epsilon^2}{\mu\mu_h},
\end{align}
which follows from the part 1 of Lemma~\ref{Lemma:Strong-Duality} and Corollary~\ref{Corollary:Strongly-convex-Lipschitz-continuous-constant-step}. Thus, claims 1 and 3 of the proposition are immediate from  \eqref{eq:Lipschitz-continuous-Srong-Convex-g-Primal-Result-1-}, together with \eqref{eq:Strongly-convexL-ipschitz-continuous-constant-step}. Moreover, we have
\begin{align}\nonumber
& f(\mathbf{y}^{(k)})-f(\mathbf{y}^\star)   
 \leq  \big(1+S\sqrt{2/\mu_h}\big)\big(h(\boldsymbol{\lambda}^{(k)})-h(\boldsymbol{\lambda}^\star)\big) \nonumber \\ \label{eq:Lipschitz-continuous-Srong-Convex-g-Primal-Result-second-2-}
 & \qquad \qquad \qquad \qquad \qquad + S\|\boldsymbol{\lambda}^{\star}\| \sqrt{h(\boldsymbol{\lambda}^{(k)})-h(\boldsymbol{\lambda}^\star)},
\end{align}
where $S=\sqrt{(4+4\cos(\pi/m))/\mu}$. Here we have used the strong convexity of $h$, \cf Proposition~\ref{Proposition:Strong-Concavity-of-g}. In particular, from \cite[p.~11, (35)]{Polyak-Intro-Opt-1987}, we have
\begin{equation}
    \|\boldsymbol{\lambda}^{(k)}-\boldsymbol{\lambda}^\star\|^2\leq \frac{2}{\mu_h} \ \left(h(\boldsymbol{\lambda}^{(k)})-h(\boldsymbol{\lambda}^\star)\right),
\end{equation}
which in turn ensures that $\|\boldsymbol{\lambda}^{(k)}\|\leq\|\boldsymbol{\lambda}^{\star}\|{+} \sqrt{(2/\mu_h)\big(h(\boldsymbol{\lambda}^{(k)})-h(\boldsymbol{\lambda}^\star)\big) }$ since $\|\boldsymbol{\lambda}^{(k)}\|-\|\boldsymbol{\lambda}^{\star}\|\leq \|\boldsymbol{\lambda}^{(k)}-\boldsymbol{\lambda}^\star\|$. Thus, claims 2 and 4 of the proposition follow from \eqref{eq:Lipschitz-continuous-Srong-Convex-g-Primal-Result-second-2-} and \eqref{eq:Strongly-convexL-ipschitz-continuous-constant-step}. 
\end{IEEEproof}

Unlike Proposition~\ref{Proposition:Lipschitz-continuous-g-Primal-Result}, in which the distance to the dual optimal solution $\|\boldsymbol{\lambda}^{(k)}-\boldsymbol{\lambda}^\star\|$ is assumed to be uniformly bounded by some $D$, in Proposition~\ref{Proposition:Strongly-convex-Lipschitz-continuous-constant-step-Primal} no such assumptions are made. This is a consequence of the strong convexity of $h$.

Convergence with nonsummable stepsize is established below, where \cite[\S~2.2, Lemma~3]{Polyak-Intro-Opt-1987} plays a key role.

\begin{corollary} \label{Corollary:Strongly-convexL-ipschitz-continuous-nonsum-step}
Suppose Assumption\ref{Assumption:Absolute-Deterministic-Error},
Assumption~\ref{Assumption:Strongly-Convex-f}, and Assumption~\ref{Assumption:Lipschitz-Convex-Grad-f} hold.  Let $\{\boldsymbol{\lambda}^{(k)}\}$ be the sequence of Lagrange multipliers generated by Algorithm~\ref{Alg:Fully-Distributed}. Moreover, suppose $\gamma_k$ satisfy the nonsummable  stepsize rule with $0<\gamma_k\leq1/L_h$. Then
\begin{enumerate}
    \item $\displaystyle \underset{k}{{\limsup}} \ h(\boldsymbol{\lambda}^{(k)})-h(\boldsymbol{\lambda}^{\star})\leq \frac{\epsilon^2}{2\mu_h}$.
    \item for $\gamma_k=(c/\mu_h)/(k+1)^{p}$, where $0<p\leq1$ and $0< c\leq\mu_h/L_h$,
\begin{equation} \label{eq:Strongly-convexL-ipschitz-continuous-nonsum-conv-rate}
h(\boldsymbol{\lambda}^{(k)})-h(\boldsymbol{\lambda}^{\star})= O\left(\frac{1}{k^{c/p}}\right) +\frac{\epsilon^2}{2\mu_h}.
\end{equation}
\end{enumerate}
\end{corollary}
%
\begin{IEEEproof}
By rearranging the terms in \eqref{eq:Strongly-convex-Lipschitz-continuous-constant-nonsum-step}, we get 
    \begin{align}
    h(&\boldsymbol{\lambda}^{(k+1)})-h(\boldsymbol{\lambda}^{\star})-{\epsilon^2}/{2\mu_h} \nonumber \\
    &\leq (1-\gamma_k \mu_h)\left(h(\boldsymbol{\lambda}^{(k)})-h(\boldsymbol{\lambda}^{\star})-{\epsilon^2}/{2\mu_h}\right). \label{eq:Strongly-convexL-ipschitz-continuous-nonsum-step-2-}
    \end{align}
    Let $q^{(k)}=1-\gamma_k \mu_h$ and $\alpha^{(k)}=0$. Thus we have $\textstyle \sum_{k=0}^\infty(1-q^{(k)})=\infty$ from the suppositions. Moreover, we have $\alpha^{(k)}/(1-q^{(k)})\to 0$ as $k\to \infty$. Hence, \cite[\S~2.2, Lemma~3]{Polyak-Intro-Opt-1987} readily applies to yield  $\limsup_{k}h(\boldsymbol{\lambda}^{(k)})-h(\boldsymbol{\lambda}^{\star})-(\epsilon^2/ 2\mu_h)\leq0$, which completes the proof of part 1.
    
    To show the second part, we let $u^{(k)}=h(\boldsymbol{\lambda}^{(k)})-h(\boldsymbol{\lambda}^{\star})-(\epsilon^2/2\mu_h)$ and  $v^{(k)}=(k+1)^{c/p} u^{(k)}$. Then,
    \begin{align}
        & v^{(k+1)}\leq \left(1+\frac{1}{k+1}\right)^{c/p}\left(1-\frac{c}{(k+1)^p}\right)v^{(k)} \allowdisplaybreaks \label{eq:Strongly-convex-and-Lipschitz-continuous-nonsum-conv-rate-3-}\\
         &\leq \left(1+\frac{c}{p(k+1)}+\frac{c^2}{2p^2(k+1)^2}+o\left(\frac{1}{(k+1)^2}\right)\right)\nonumber\\
         &\qquad \qquad \qquad \qquad \qquad  \qquad \times\left(1-\frac{c}{(k+1)^p}\right)v^{(k)}\allowdisplaybreaks \label{eq:Strongly-convex-and-Lipschitz-continuous-nonsum-conv-rate-4-}\\
          &= \left(1-\frac{c^2}{2(k+1)^{2p}}+o\left(\frac{1}{(k+1)^2}\right)\right)v^{(k)},\ \text{for all} \nonumber \\
          & \qquad \qquad \qquad \qquad \qquad  \qquad k\geq \lceil e^{\frac{\log p}{(p-1)}}-1\rceil  \allowdisplaybreaks \label{eq:Strongly-convex-and-Lipschitz-continuous-nonsum-conv-rate-8-}\\
          &\leq v^{(k)},\ \text{for sufficiently large} \ k, \label{eq:Strongly-convex-and-Lipschitz-continuous-nonsum-conv-rate-9-}
    \end{align}
where \eqref{eq:Strongly-convex-and-Lipschitz-continuous-nonsum-conv-rate-3-} follows from \eqref{eq:Strongly-convexL-ipschitz-continuous-nonsum-step-2-}, the definition of $v^{(k)}$, the choice of $\gamma_k$, together with some algebraic manipulations. The inequality \eqref{eq:Strongly-convex-and-Lipschitz-continuous-nonsum-conv-rate-4-} follows from binomial expansion and \eqref{eq:Strongly-convex-and-Lipschitz-continuous-nonsum-conv-rate-8-} is immediate from that $(x+1)^p\leq p(x+1)$ for sufficiently large $x$, when $p\in(0,1]$. Summing  \eqref{eq:Strongly-convex-and-Lipschitz-continuous-nonsum-conv-rate-9-} over $k$, together with the definitions of $v^{(k)}$ and $u^{(k)}$ yields the result.
\end{IEEEproof}

Corollary~\ref{Corollary:Strongly-convexL-ipschitz-continuous-nonsum-step} indicates that the least upperbound of $h(\boldsymbol{\lambda}^{(k)})$ converges into a neighborhood of the optimal value $h(\boldsymbol{\lambda}^{\star})$ at a rate of $O(1/k^{c/p})$, where the size of the neighborhood explicitly depends on $\epsilon$ and the strong convexity constant $\mu_h$ of $h$. Note that the rate of convergence depends on the ratio $c/p$. {It can easily be observed that, for a given $c$ value, the rate of convergence increases when the value of $p$ decreases, where $0< c\leq\mu_h/L_h$ and $0<p\leq1$.} 

Note that the case $p=0$ corresponds to a constant stepsize rule. This suggests, as in the constant stepsize rule [\cf Corollary~\ref{Corollary:Strongly-convex-Lipschitz-continuous-constant-step}], that when $p\to 0$, \eqref{eq:Strongly-convexL-ipschitz-continuous-nonsum-conv-rate} should be a good resemblance of \eqref{eq:Strongly-convexL-ipschitz-continuous-constant-step}. According to the proof of Corollary~\ref{Corollary:Strongly-convexL-ipschitz-continuous-nonsum-step}, $p$ is to be chosen in such a manner that $(k+1)^p\leq p(k+1)$ for sufficiently large~$k$. One such choice is $p=\log k/k$. The following remark will shed some light on this possibility.

\begin{remark}\label{Remark:Corollary:Strongly-convexL-ipschitz-continuous-nonsum-step-p-is-zero}
Let $s(k;p)\in O(1/k^{c/p})$ where $c>0$ and $0<p<1$. Then $s(k;p)$ with $p=\log k/k$ converges to $0$ with the rate of geometric progression.
\end{remark}
\begin{IEEEproof}
This is clear from the notation for big $O$.
\end{IEEEproof}
%
%
\begin{remark}\label{Remark:Strongly-convex-Lipschitz-continuous-nonsummable-step-Primal}
By using Corollary~\ref{Corollary:Strongly-convexL-ipschitz-continuous-nonsum-step}, convergence assertions similar to Proposition~\ref{Proposition:Strongly-convex-Lipschitz-continuous-constant-step-Primal} for the sequences $\{\mathbf{y}^{(k)}\}$ and $\{f(\mathbf{y}^{(k)})\}$ can be derived analogously for nonsummable $\gamma_k$. 
\end{remark}

\subsection{Feasible Points from Algorithm~\ref{Alg:Fully-Distributed}}\label{subsec:Feasible-Points-and-Convergence}

Proposition~\ref{Proposition:Lipschitz-continuous-g-Primal-Result} and~\ref{Proposition:Strongly-convex-Lipschitz-continuous-constant-step-Primal} characterize how far from primal solution $\mathbf{y}^\star$ the locally computed solution $\mathbf{y}^{(k)}$ is. However, $\mathbf{y}^{(k)}$ is not necessarily feasible for problem~\eqref{eq:distributed-problem}, despite~$k$ being very large. More specifically, $\mathbf{A}\mathbf{y}^{(k)}\neq \mathbf{0}$. From both analytical and practical standpoints, the computation of a \emph{feasible point}, together with its convergence properties is also of crucial importance. To this end, the following remark will be useful. 
\begin{remark}\label{Remark:Compute-Feasible-Point}
Let ${\tilde{\mathbf{y}}}^{(k)}$ be a point in $\R^{nm}$ given by
\begin{equation}
    {\tilde {\mathbf{y}}}^{(k)} = ({1}/{m}) \ (\boldsymbol{1}_{m\times m}\otimes \mathbf{I}_n) \ {\mathbf{y}}^{(k)},
\end{equation}
where $\mathbf{y}^{(k)}$ is given in \eqref{eq:Local-y-k-in-vector-form}. Then ${\tilde{ \mathbf{y}}}^{(k)}$ is a feasible point for problem~\eqref{eq:distributed-problem}.
\end{remark}
\begin{IEEEproof}
The proof is straightforward and is omitted.
\end{IEEEproof}

\begin{lemma}\label{Lemma:Feasible-Points}
Let $\mathbf{y}$ be a vector in $\bar{\mathcal{Y}}$ given in \eqref{eq:Cartisian-Profuct-Y} and $\tilde{\mathbf{y}} = {1}/{m} \ (\boldsymbol{1}_{m\times m}\otimes \mathbf{I}_n) \ \mathbf{y}$. Moreover, suppose the optimal value $p^\star$ of problems \eqref{eq:distributed-problem} is attained for some $\mathbf{y}^\star\in\{\mathbf{y}\in\R^{nm} \ | \ \forall~i \  \mathbf{y}_i\in\mathcal{Y}, \  \mathbf{A}\mathbf{y}=\mathbf{0} \}$. Then, 
\begin{enumerate}
    \item   $ \| \mathbf{y}-\mathbf{y}^\star\|\geq \| \tilde{ \mathbf{y}}-\mathbf{y}^\star\|$. 
\item $ \tilde D\| \tilde {\mathbf{y}}-\mathbf{y}^\star\|\geq f(\tilde {\mathbf{y}})-f(\mathbf{y}^\star)$, where
\end{enumerate}
\begin{equation}\label{eq:Primal-Subgrad-Uniform-Bound}
    \tilde D = {\sup}_{\mathbf{A}\hat {\mathbf{y}}=\mathbf{0},\hat {\mathbf{y}}\in\bar{\mathcal{Y}}} \ \{\|\boldsymbol{\nu}\| \ | \ \boldsymbol{\nu}\in\partial f(\hat {\mathbf{y}}) \},
\end{equation}
and $\mathbf{A}$ is defined in \eqref{eq:A-Matrix-for-Consensus-Constraint}.
\end{lemma}
\begin{IEEEproof}
See Appendix \ref{Appendix:Lemma:Feasible-Points}.
\end{IEEEproof}
Note that part 2 of Lemma~\ref{Lemma:Feasible-Points} relies on certain Lipschitzian properties of the primal function~$f$, \cf Assumption \ref{Assumption:Lipschitz-Convex-Grad-f}. 
Let us next establish convergences of $\{\tilde{\mathbf{y}}^{(k)}\}$.
\begin{prop2}\label{Proposition:Feasible-Point-Lipschitz-continuous-g-Primal-Result}
Suppose Assumption~\ref{Assumption:Absolute-Deterministic-Error}, Assumption~\ref{Assumption:Strong-Duality}, and Assumption~\ref{Assumption:Strongly-Convex-f} hold. Let $\{\boldsymbol{\lambda}^{(k)}\}$ be the sequence of Lagrange multipliers generated by Algorithm~\ref{Alg:Fully-Distributed}, $\{\mathbf{y}^{(k)}\}$ be the corresponding sequence of primal variables, and $\{\tilde {\mathbf{y}}^{(k)}\}$ be the resulting sequence of primal feasible points, where $\tilde {\mathbf{y}}^{(k)} = ({1}/{m}) \ (\boldsymbol{1}_{m\times m}\otimes \mathbf{I}_n) \ \mathbf{y}^{(k)}$. Moreover, suppose that the functions $f_i$, $i=1,\ldots,m$ are differentiable and the distance to the dual optimal solution $\|\boldsymbol{\lambda}^{(k)}-\boldsymbol{\lambda}^\star\|$ is uniformly bounded by some $D$. Let $\gamma_k$ satisfy the nonsummable  stepsize rule with $0<\gamma_k\leq1/L_h$. Then
\begin{enumerate}
    \item $\displaystyle \underset{k}{{\limsup}} \ \underset{i\in\{0,\ldots,k\}}{\min} \ \|\tilde {\mathbf{y}}^{(i)}-\mathbf{y}^\star\| \leq \rho_{1}(\epsilon)$, where $\rho_{1}(\epsilon)\to 0$ as $\epsilon\to 0$.
    \item for $\gamma_k=\gamma/(k+1)^p$, where $0<\gamma\leq1/L_h$ and $0\leq p\leq1$,
\begin{equation*}\label{eq:Lipschitz-continuous-nonsum-Step-convergence-rate}
\underset{i\in\{0,\ldots,k\}}{\min} \ \|\tilde {\mathbf{y}}^{(i)}-\mathbf{y}^\star\| =
    \begin{cases}
    \displaystyle O\left(\frac{1}{\sqrt[4]{k^{{1-p}}}}\right)+\rho_{1}(\epsilon),& \\
    & \hspace{-7.5mm} \text{$p\in[0,1)$}\\
         \displaystyle O\left( \frac{1}{\sqrt[4]{\log k}} \right) +\rho_{1}(\epsilon), & \\
         & \hspace{-3.5mm} \text{$p=1$},
    \end{cases}       
\end{equation*}
and the best convergence rate is of the order $O(1/\sqrt[4]{k})$, which is achieved when $p=0$.
\end{enumerate}
\end{prop2}
\begin{IEEEproof}
Combining  parts 1 and 3 of Proposition~\ref{Proposition:Lipschitz-continuous-g-Primal-Result}, together with part 1 of Lemma~\ref{Lemma:Feasible-Points}, yields the result.
\end{IEEEproof}

We did not write explicitly the underlying expression for $\rho_1$ for brevity. The end results would essentially be the same as those claimed in Proposition~\ref{Proposition:Lipschitz-continuous-g-Primal-Result}.
We note that the hypotheses of Proposition~\ref{Proposition:Feasible-Point-Lipschitz-continuous-g-Primal-Result} do not permit any means for quantifying a bound on the error of primal objective evaluated at feasible points. This restriction could not be avoided due to technical reasons. More specifically, the least upperbound on the norm of the subgradients $\boldsymbol\nu$ of $f$ evaluated over the set $\{\hat {\mathbf{y}} \ | \ \mathbf{A}\hat {\mathbf{y}}=\mathbf{0},\hat {\mathbf{y}}\in\bar{\mathcal{Y}}\}$ [\cf \eqref{eq:Primal-Subgrad-Uniform-Bound}] can be $\infty$ because the function $f$ might have a nontrivial domain. A quantification is possible nevertheless, as will be asserted in the next proposition, if in addition, the gradient Lipschitz property of $f$ stated in Assumption~\ref{Assumption:Lipschitz-Convex-Grad-f} is imposed.

\begin{prop2}\label{Proposition:Feasible-Point-convex-Lipschitz-continuous-constant-step-Primal}
Suppose Assumption~\ref{Assumption:Absolute-Deterministic-Error}, Assumption~\ref{Assumption:Strong-Duality}, Assumption~\ref{Assumption:Strongly-Convex-f}, and Assumption~\ref{Assumption:Lipschitz-Convex-Grad-f} hold. Let $\{\boldsymbol{\lambda}^{(k)}\}$ be the sequence of Lagrange multipliers generated by Algorithm~\ref{Alg:Fully-Distributed}, $\{\mathbf{y}^{(k)}\}$ be the corresponding sequence of primal variables, and $\{\tilde {\mathbf{y}}^{(k)}\}$ be the resulting sequence of primal feasible points, where $\tilde {\mathbf{y}}^{(k)} = ({1}/{m}) \ (\boldsymbol{1}_{m\times m}\otimes \mathbf{I}_n) \ \mathbf{y}^{(k)}$. Let the stepsize $\gamma_k=\gamma$ for all $k\in\mathbb{Z}^0_+$. 
Then for $0<\gamma\leq1/L_h$
\begin{enumerate}
    \item $\displaystyle \underset{k}{{\limsup}} \ \|\tilde {\mathbf{y}}^{(k)}-\mathbf{y}^\star\| \leq \rho_{2}(\epsilon)$, where $\rho_{2}(\epsilon)\to 0$ as $\epsilon\to 0$.
\item  $\displaystyle \underset{k}{{\limsup}} \ \|f(\tilde {\mathbf{y}}^{(k)})-f(\mathbf{y}^\star)\| \leq  \rho_{3}(\epsilon)$, where $\rho_{3}(\epsilon)\to 0$ as $\epsilon\to 0$.
    \item  the least upperbound of $ \|\tilde {\mathbf{y}}^{(k)}-\mathbf{y}^\star\|$ converges into a neighborhood of $0$ with the rate of geometric progression. 
 \item the least upperbound of $ \|f(\tilde {\mathbf{y}}^{(k)})-f(\mathbf{y}^\star)\|$ converges into a neighborhood of $0$ with the rate of geometric progression.
\end{enumerate}
\end{prop2}
\begin{IEEEproof}
The first and the second claims follow from parts 1-2 of Proposition~\ref{Proposition:Strongly-convex-Lipschitz-continuous-constant-step-Primal} and parts 1-2 of Lemma~\ref{Lemma:Feasible-Points}. The proof of the last two claims is identical to that of the preceding Proposition~\ref{Proposition:Strongly-convex-Lipschitz-continuous-constant-step-Primal} and thus omitted.
\end{IEEEproof}

\begin{remark}\label{Remark:Feasible-Point-convex-Lipschitz-continuous-nonsumable-step-Primal}
By using the points highlighted in Corollary~\ref{Corollary:Strongly-convexL-ipschitz-continuous-nonsum-step} and Remark~\ref{Remark:Strongly-convex-Lipschitz-continuous-nonsummable-step-Primal}, 
convergence assertions similar to Proposition~\ref{Proposition:Feasible-Point-convex-Lipschitz-continuous-constant-step-Primal} for the sequences of feasible points $\{\tilde {\mathbf{y}}^{(k)}\}$ and 
feasible objective values $\{f(\tilde {\mathbf{y}}^{(k)})\}$ can be derived analogously for nonsummable $\gamma_k$. 
\end{remark}

\section{Numerical Results} \label{sec:Numerical-Results}

We next verify empirically the theoretical assertions in \S~\ref{sec:Convergence-Analysis-I}. To this end, problem \eqref{eq:distributed-problem} is considered with quadratic $f_i$s, i.e.,
\begin{equation} \label{eq:Numerical-section-local-problems}
    f_i(\mathbf{y}_i)= \mathbf{y}_i\tran \mathbf{A}_i \mathbf{y}_i+\mathbf{q}_i\tran \mathbf{y}_i, \quad \mathbf{A}_i\in \mathbb{S}^{n}_{++}, \ \mathbf{q}_i\in\R^n,
\end{equation}
where $\mathbb{S}^{n}_{++}$ is the \emph{positive definite cone}, and $\mathbf{A}_i$ and $\mathbf{q}_i$ are arbitrarily chosen. The Lipshitz and strong convexity constants associated with $h=-g$ (i.e., $L_h$ and $\mu_h$) are computed accordingly. Thus, Assumptions~\ref{Assumption:Strongly-Convex-f}, \ref{Assumption:Strong-Duality}, and closedness of $f_i$s hold throughout the rest of the section. Moreover, the considered system is with $n=1$ and $m=5$.


\subsection{Dual Function $g$ with Lipschitz Continuous Gradients}

Let $\mathcal{Y}=\{\mathbf{u}\in\R^n\ \, | \, -a\boldsymbol{1}_{n\times 1}\leq \mathbf{u}\leq a\boldsymbol{1}_{n\times 1}\}$, where $a>0$, i.e., $\mathcal{Y}$ is not only closed, but also compact. As a result, the dual function $g$ is with Lipschitz continuous gradient, \cf~\eqref{eq:CASE-1}. 

We consider that the distorted vector $\hat{\mathbf{d}}^{(k)}$ [\cf \eqref{eq:Distorted-Subgradient}] is a consequence of a naive quantization scheme implemented in step 4 of Algorithm~	\ref{Alg:Fully-Distributed}. In particular, the box $\mathcal{Y}$ is partitioned into identical mini-boxes of width $t=2a/2^b$ per-dimension, where $b\in\mathbb{Z}_+$. The indexing of the mini-boxes is common to all SSs. At step 4 of the algorithm, SS~$i$ first chooses $\hat {\mathbf{y}}_i^{(k)}$ to be the centroid of the mini-box in which $\mathbf{y}_i^{(k)}$ lies. Then the index of the chosen mini-box is transmitted which is simply an $nb$-bit word. As a result, the distortion $\mathbf{r}_i^{(k)}$ of $\mathbf{y}_i^{(k)}$ is bounded, i.e., $\|\mathbf{r}_i^{(k)}\|\leq \varepsilon_i=\sqrt{n}\,t/2$, conforming to Assumption~\ref{Assumption:Absolute-Deterministic-Error}. The overall norm distortion is bounded as stated below.

\begin{figure}[!t]
\centering
\subfloat[]{\includegraphics[width=0.5\linewidth]{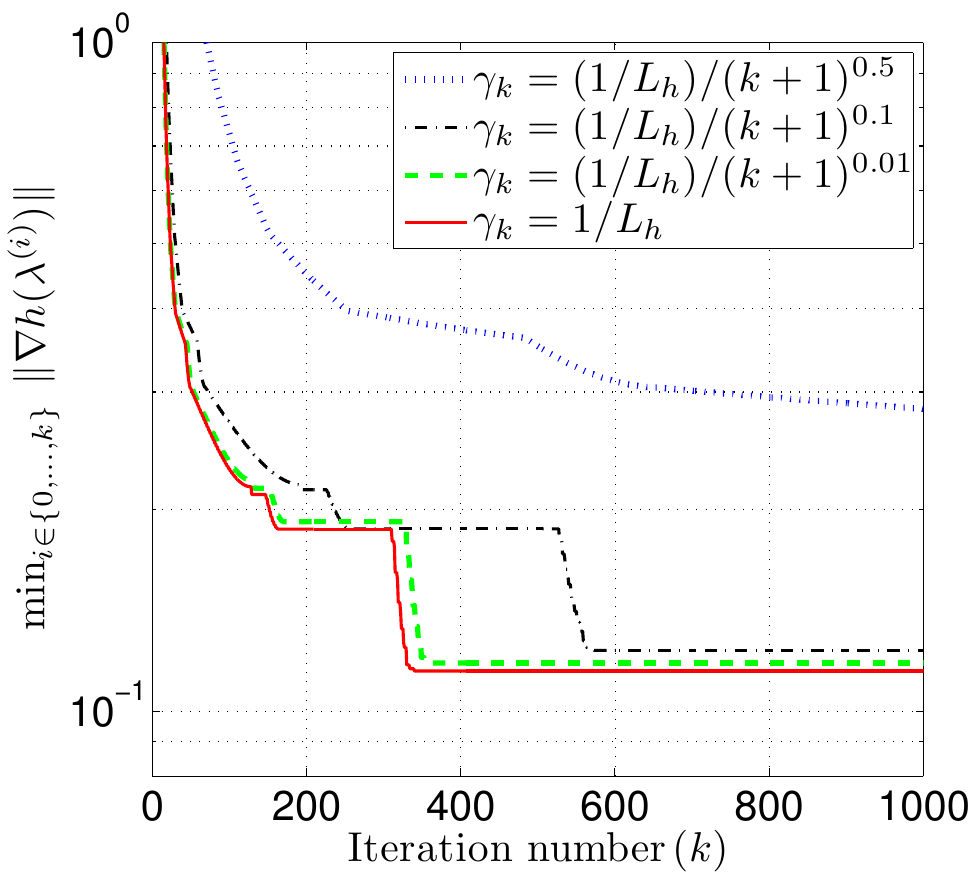}%
\label{Fig:Case1-dual-rate-3}}
\hfil
\subfloat[]{\includegraphics[width=0.5\linewidth]{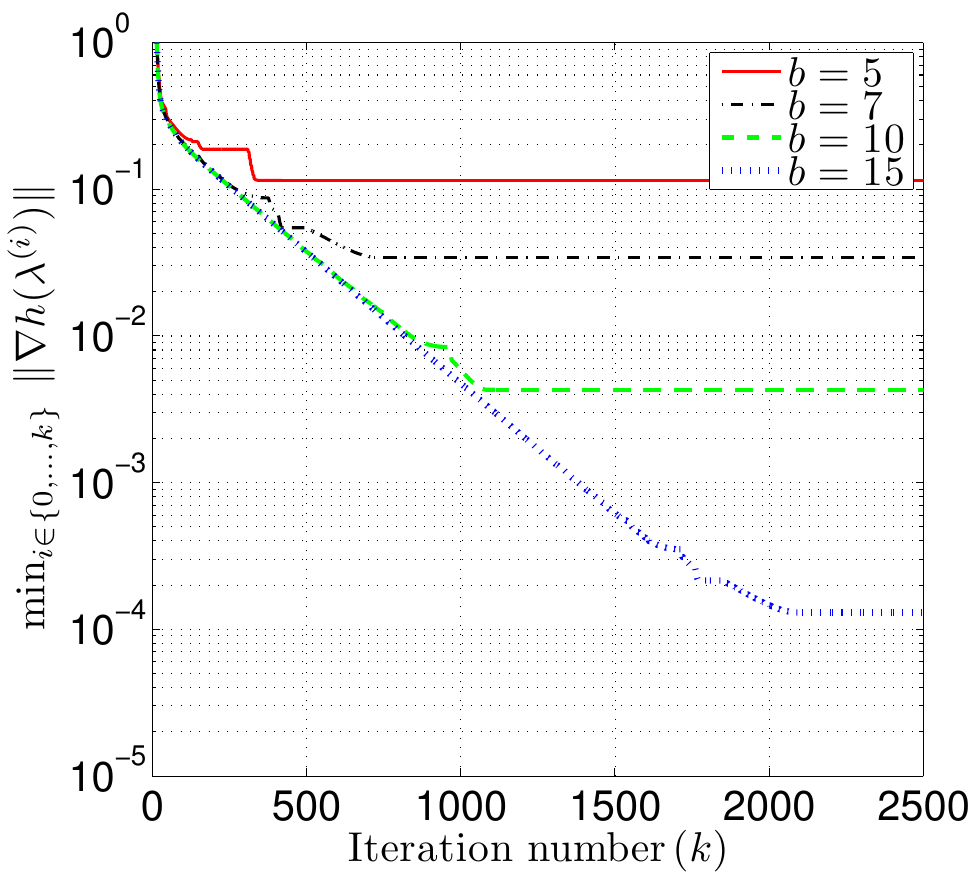}%
\label{Fig:Case1-dual-ball-4}}
\\
\subfloat[]{\includegraphics[width=0.5\linewidth]{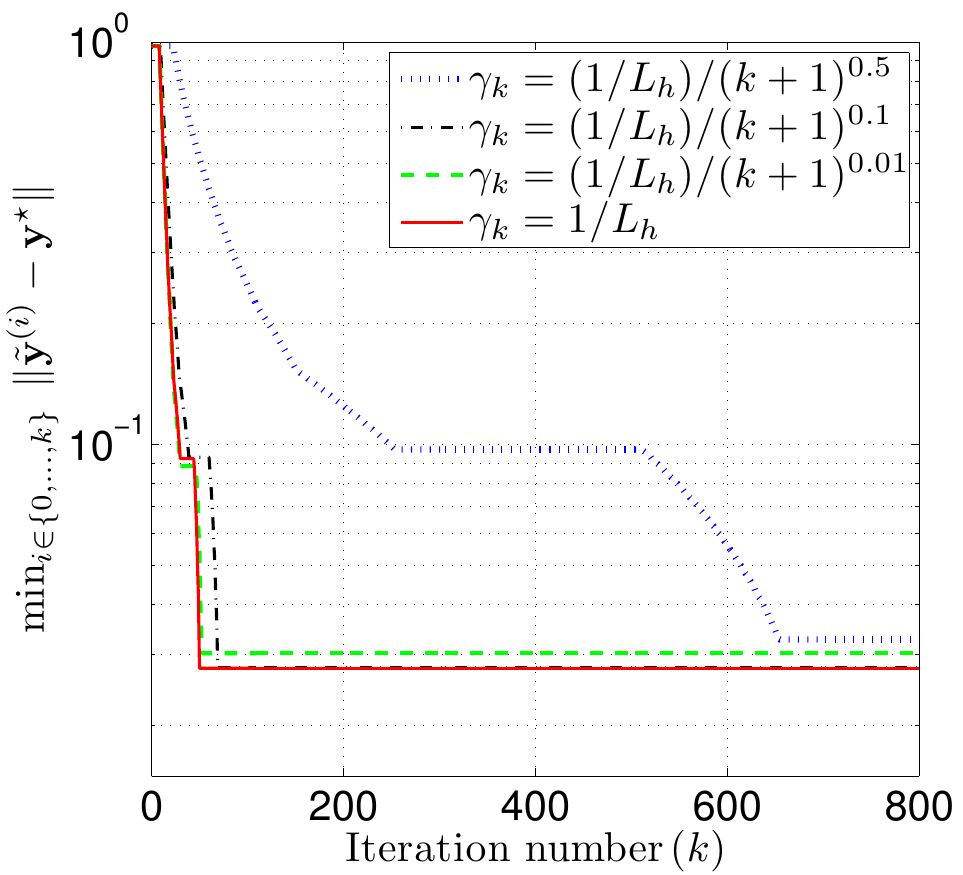}%
\label{Fig:Case1-primal-rate-1}}
\hfil
\subfloat[]{\includegraphics[width=0.5\linewidth]{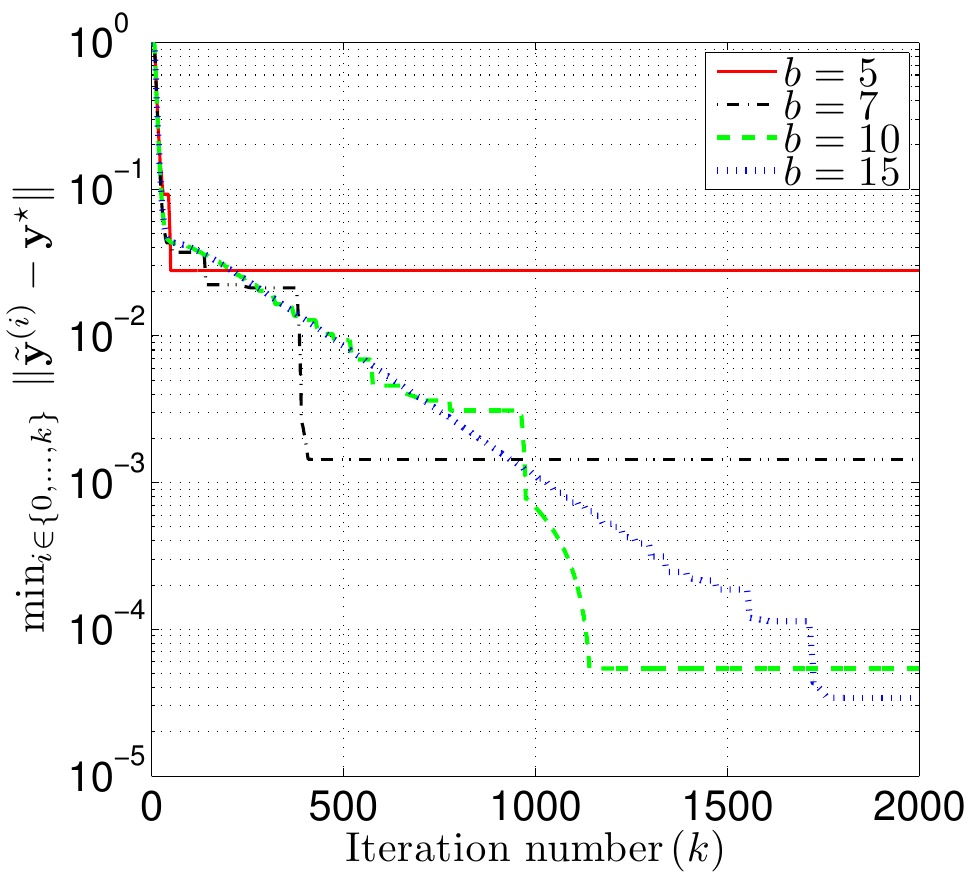}%
\label{Fig:Case1-primal-ball-2}}
\caption{\textsc{Case 1}: Convergence of dual gradients and primal feasible points, where $a=3$ in $\mathcal{Y}$ definition. \protect\subref{Fig:Case1-dual-rate-3} Effect of choice of $p$ on the convergence of dual gradients. \protect\subref{Fig:Case1-dual-ball-4} Effect of choice of $b$ on the convergence of dual gradients. \protect\subref{Fig:Case1-primal-rate-1} Effect of choice of $p$ on the convergence of primal feasible points. \protect\subref{Fig:Case1-primal-ball-2} Effect of choice of $b$ on the convergence of primal feasible points.}  
\label{Figure:Case1}
\end{figure}

\begin{remark} \label{Remark:case1-numerical-results-epsilon}
The norm of the total error vector $\mathbf{r}^{(k)}$ is bounded, i.e., $\|\mathbf{r}^{(k)}\|\leq \epsilon=6\sqrt{n(m-1)}/2^b$, \cf Corollary~\ref{Corollary:Lipschitz-continuous-nonsum-step}. 
\end{remark}

\figurename~\ref{Figure:Case1}\subref{Fig:Case1-dual-rate-3}-\subref{Fig:Case1-dual-ball-4} show the convergence of gradients of the dual function (i.e., Corollary~\ref{Corollary:Lipschitz-continuous-nonsum-step}). \figurename~\ref{Figure:Case1}\subref{Fig:Case1-dual-rate-3} depicts the effect of the choice of $p$ in the stepsize $\gamma_k=(1/L_h)/(k+1)^p$ by fixing $b=5$.
Results show that the smaller the value of $p$, the higher the rate of convergence, as claimed in Corollary~\ref{Corollary:Lipschitz-continuous-nonsum-step}-(2). Moreover, the best rate is achieved when $p=0$, which corresponds to the fixed stepsize rule. \figurename~\ref{Figure:Case1}\subref{Fig:Case1-dual-ball-4} shows the effect of the choice of $b$ by fixing $\gamma_k=1/L_h$, the fixed stepsize rule. Results show that when the number of bits $b$ increases the size of the neighborhood around $0$ to which $\min_{i}\|\nabla h(\boldsymbol{\lambda}^{(i)})\|$ converges decreases. This is readily expected from Corollary~\ref{Corollary:Lipschitz-continuous-nonsum-step}-(1), together with Remark~\ref{Remark:case1-numerical-results-epsilon}, because $\epsilon$, the neighborhood, is inversely proportional to $2^b$. 
\figurename~\ref{Figure:Case1}\subref{Fig:Case1-primal-rate-1}-\subref{Fig:Case1-primal-ball-2} show the convergence of corresponding primal feasible points (i.e., Proposition~\ref{Proposition:Feasible-Point-Lipschitz-continuous-g-Primal-Result}). Results are similar to that of \figurename~\ref{Figure:Case1}\subref{Fig:Case1-dual-rate-3}-\subref{Fig:Case1-dual-ball-4}. \figurename~\ref{Figure:Case1}\subref{Fig:Case1-primal-rate-1} shows that the number of iterations required for reaching the neighborhood in the primal-domain appears to be relatively smaller than that in the dual-domain, especially for smaller $p$ values. This behavior is typical for many methods in general, because a good feasible point can usually be computed, even with a relatively smart heuristic method.


\subsection{Strongly Concave $g$ with Lipschitz Continuous Gradients} \label{subsec:Strongly-convex-Lipschitz-continuous}

Let $\mathcal{Y}=\R^n$, i.e., $\mathcal{Y}$ is closed. In this setting, it can be verified that the dual function $g$ is strongly concave and with Lipschitz continuous gradient, \cf~\eqref{eq:CASE-2}.

We consider that the distorted vector $\hat{\mathbf{d}}^{(k)}$ [\cf \eqref{eq:Distorted-Subgradient}] is a consequence of measurement errors at CN in the partially distributed variant, \cf Figure~\ref{Fig:Partially-Distributed-Algorithm}. The magnitudes of measurement errors are bounded from above by some $\varsigma>0$ per dimension. As a result, the distortion $\mathbf{r}_i^{(k)}$ of $\mathbf{y}_i^{(k)}$ is bounded, i.e., $\|\mathbf{r}_i^{(k)}\|\leq \varepsilon_i=\sqrt{n}\,\varsigma$, conforming to Assumption~\ref{Assumption:Absolute-Deterministic-Error}. The norm distortion is bounded as follows.
\begin{remark} \label{Remark:case2-numerical-results-epsilon}
The norm of the total error vector $\mathbf{r}^{(k)}$ is bounded, i.e., $\|\mathbf{r}^{(k)}\|\leq \epsilon=2\sqrt{n(m-1)}\, \varsigma$, \cf Corollary~\ref{Corollary:Strongly-convex-Lipschitz-continuous-constant-step} and Corollary~\ref{Corollary:Strongly-convexL-ipschitz-continuous-nonsum-step}. 
\end{remark}

\begin{figure}[!t]
\centering
\subfloat[]{\includegraphics[width=0.5\linewidth]{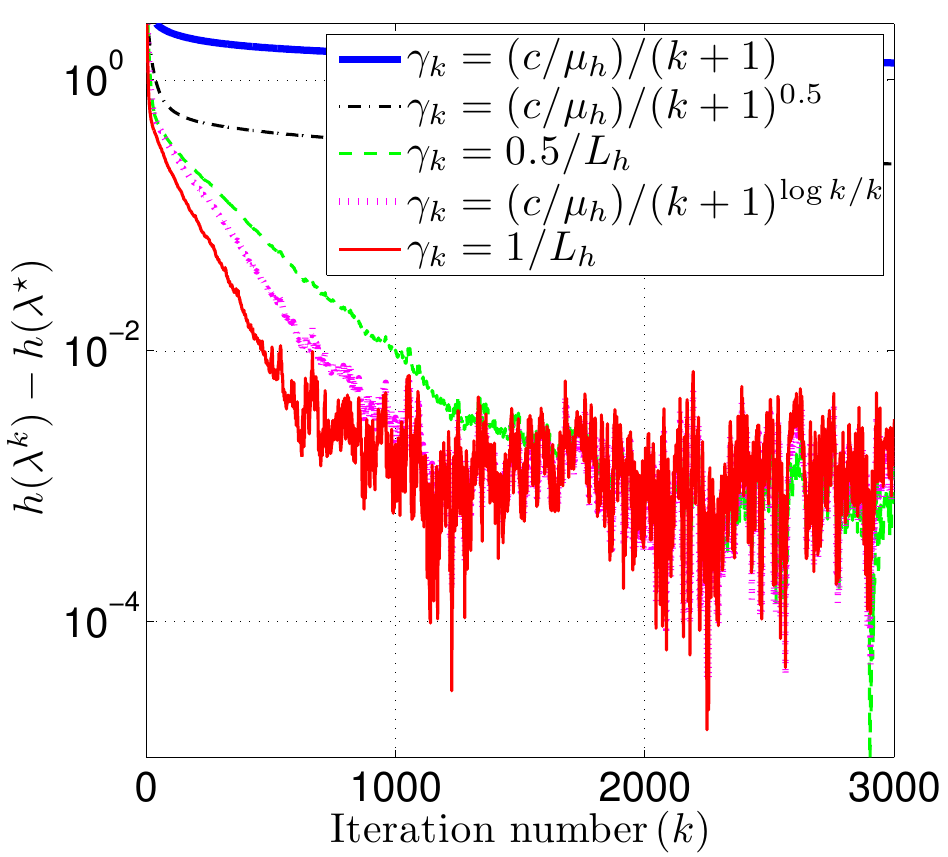}%
\label{Fig:Case2-dual-rate-3}}
\hfil
\subfloat[]{\includegraphics[width=0.5\linewidth]{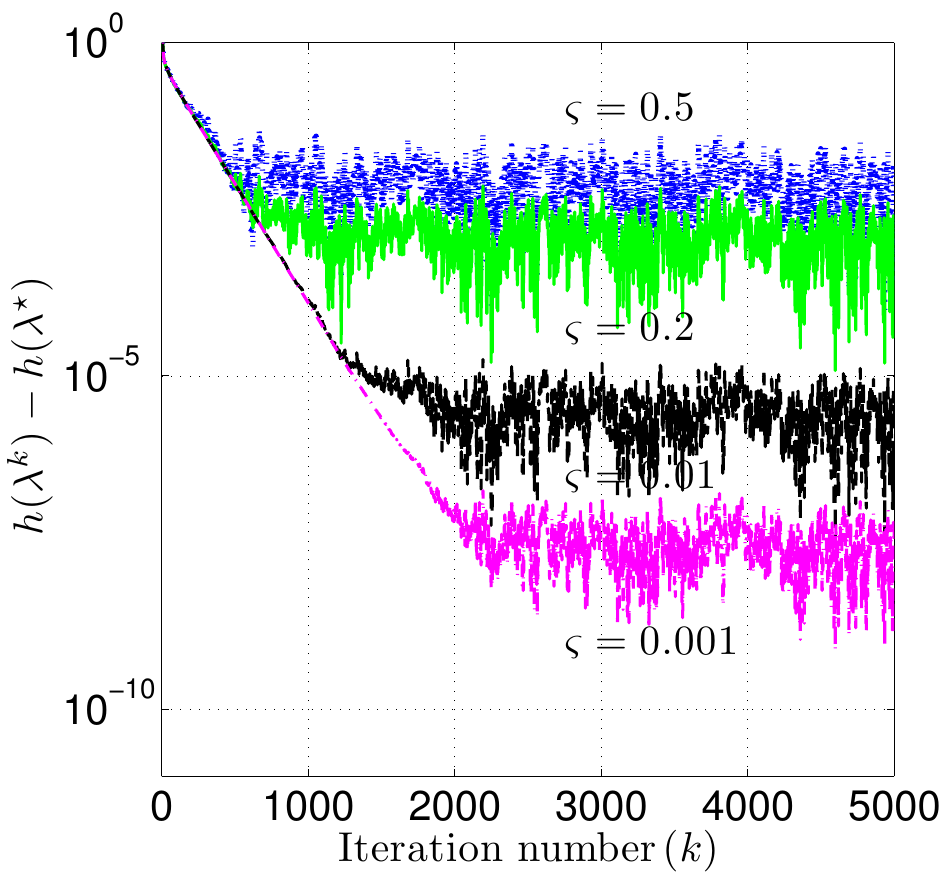}%
\label{Fig:Case2-dual-ball-4}}
\\
\subfloat[]{\includegraphics[width=0.5\linewidth]{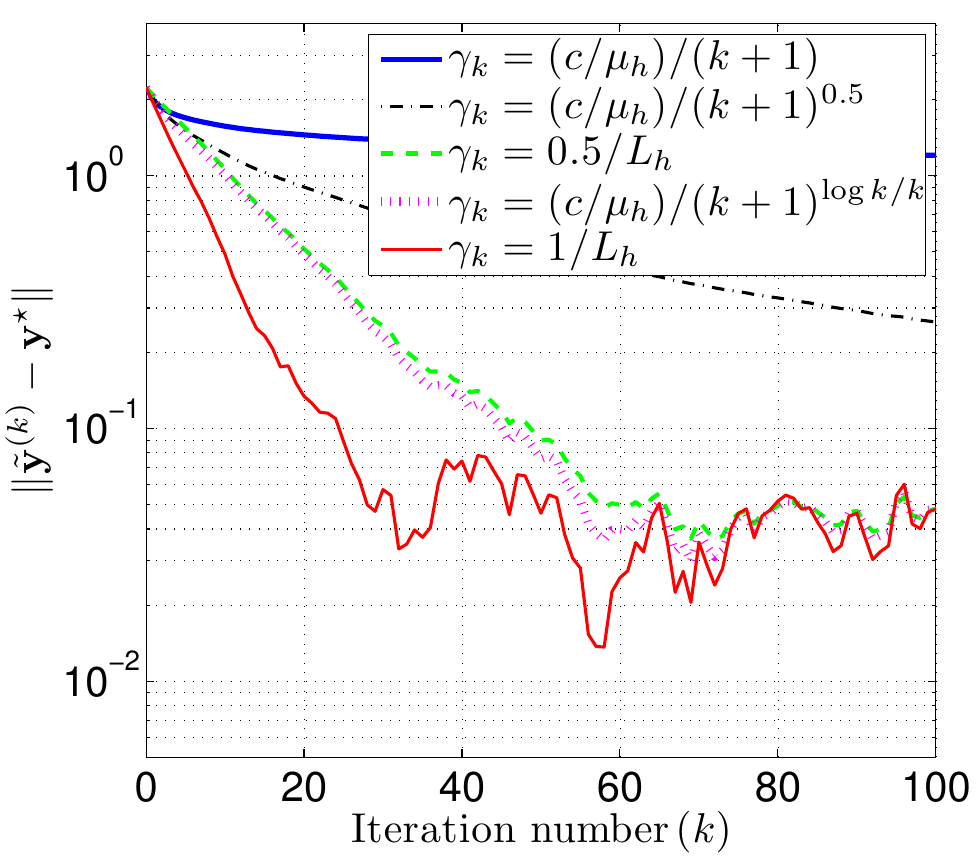}%
\label{Fig:Case2-primal-rate-1}}
\hfil
\subfloat[]{\includegraphics[width=0.5\linewidth]{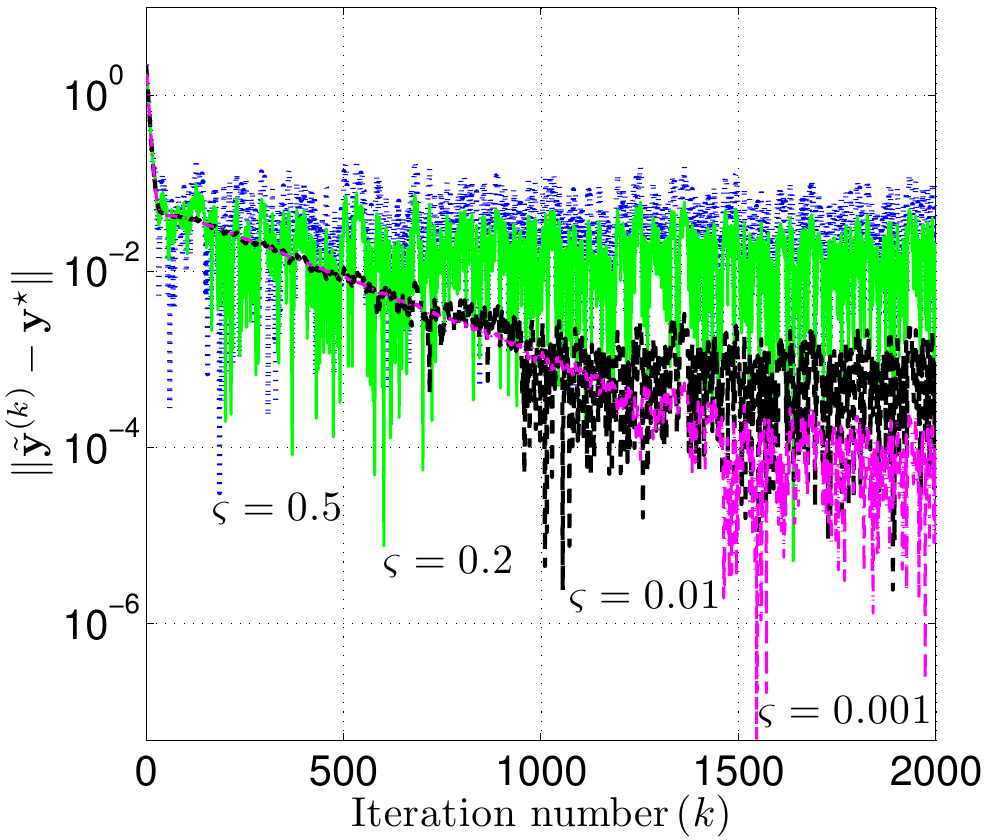}%
\label{Fig:Case2-primal-ball-2}}
\caption{\textsc{Case 2}: Convergence of dual function values and primal feasible points. \protect\subref{Fig:Case2-dual-rate-3} Dual function values using constant and nonsummable stepsizes. \protect\subref{Fig:Case2-dual-ball-4} Effect of choice of $\varsigma$ on the convergence of dual function  values. \protect\subref{Fig:Case2-primal-rate-1} Primal feasible points using constant and nonsummable stepsizes. \protect\subref{Fig:Case2-primal-ball-2} Effect of choice of $\varsigma$ on the convergence of primal feasible points.} 
\label{Figure:Case2}
\end{figure}

\figurename~\ref{Figure:Case2}\subref{Fig:Case2-dual-rate-3}-\subref{Fig:Case2-dual-ball-4} show the convergence of dual function values for both the fixed stepsize rule (i.e., Corollary~\ref{Corollary:Strongly-convex-Lipschitz-continuous-constant-step}) and nonsummable stepsize rule (i.e., Corollary~\ref{Corollary:Strongly-convexL-ipschitz-continuous-nonsum-step}). \figurename~\ref{Figure:Case2}\subref{Fig:Case2-dual-rate-3} shows that linear convergence is guaranteed with fixed stepsizes, while $\gamma_k=1/L_h$ being the best choice. This clearly agrees with the assertions claimed in Corollary~\ref{Corollary:Strongly-convex-Lipschitz-continuous-constant-step}, \cf~\eqref{eq:Strongly-convexL-ipschitz-continuous-constant-step}. Moreover, the results demonstrate the effect of the choice of $p$ in the stepsize $\gamma_k=(c/\mu_h)/(k+1)^p$ by fixing $\varsigma=0.2$ and $c=0.004$. Note that $c$ is carefully chosen so that it lies inside the prescribed limits $0< c\leq\mu_h/L_h$ imposed by Corollary~\ref{Corollary:Strongly-convexL-ipschitz-continuous-nonsum-step}. Results show that the smaller the value of $p$, the higher the rate of convergence, as claimed in Corollary~\ref{Corollary:Strongly-convexL-ipschitz-continuous-nonsum-step}-(2). For comparisons, we have also included the convergence of dual function values for $\gamma_k=(c/\mu_h)/(k+1)^p$ with $p=\log k/k$, which can be interpreted as a limiting case of $\gamma_k=(c/\mu_h)/(k+1)^p$ as $p\to 0$, \cf Remark~\ref{Remark:Corollary:Strongly-convexL-ipschitz-continuous-nonsum-step-p-is-zero}. Results clearly demonstrate a linear convergence, as claimed in Remark~\ref{Remark:Corollary:Strongly-convexL-ipschitz-continuous-nonsum-step-p-is-zero}. 
\figurename~\ref{Figure:Case2}\subref{Fig:Case2-dual-ball-4} shows the effect of the choice of $\varsigma$ by fixing $\gamma_k=1/L_h$. Results show that when $\varsigma$ decreases, so is the size of the neighborhood around $h(\boldsymbol{\lambda}^\star)$ to which $h(\boldsymbol{\lambda}^{(k)})$ converges. This behavior is expected from Corollary~\ref{Corollary:Strongly-convex-Lipschitz-continuous-constant-step}, together with Remark~\ref{Remark:case2-numerical-results-epsilon}, because $\epsilon$ that defines the neighborhood, is linearly related to $\varsigma$. 
\figurename~\ref{Figure:Case2}\subref{Fig:Case2-primal-rate-1}-\subref{Fig:Case2-primal-ball-2} show the convergence of corresponding primal feasible points (i.e., Proposition~\ref{Proposition:Feasible-Point-convex-Lipschitz-continuous-constant-step-Primal}, Remark~\ref{Remark:Feasible-Point-convex-Lipschitz-continuous-nonsumable-step-Primal}). With respect to the rate of convergence [\cf \figurename~\ref{Figure:Case2}\subref{Fig:Case2-primal-rate-1}] and the size of the converging neighborhood [\cf \figurename~\ref{Figure:Case2}\subref{Fig:Case2-primal-ball-2}], results demonstrate a similar behavior to that of  \figurename~\ref{Figure:Case2}\subref{Fig:Case2-dual-rate-3}-\subref{Fig:Case2-dual-ball-4}.


\section{Conclusion} \label{Conclusion}

The \emph{global consensus} optimization problem, which plays a key role in many large-scale signal processing and machine learning application domains was considered. An inexact fully distributed algorithm is provided where the \emph{inexactness} is considered to be \emph{additive and bounded}, which in turn modeled a wide range of distortions, including quantization errors, approximation errors, errors due to subproblem solver accuracy, errors in dual variable coordination, noise in wireless settings, and measurement errors, among others. Convergences of dual and primal \emph{feasible points}, together with their \emph{rates of convergences}, were extensively analyzed. Our analytical assertions showed that the feasible points converge into a \emph{neighborhood} of optimality. The \emph{size of the neighborhood} was explicitly quantified in terms of the underlying inexactness. 
Numerical experiments were conducted to verify the theoretical assertions. Future research will include extending the results to nondifferentiable settings. Furthermore, it would be more interesting to seek how the implications of the results, with extensions, if any, can be used when designing and analyzing wireless networks for machine learning, which is of crucial importance.

\begin{appendices}

\section{Lipschitzian, Strong Convexity, and Other Properties of the Dual Function}\label{sec:Characteristics-of-Dual-Function}

We start by highlighting a useful relationship between $g$ [\cf~\eqref{eq:dual-function}], and the conjugate function~\cite[p. 473]{Rockafellar-98} of $f+\delta_{\bar{\mathcal{Y}}}$~\footnote{Recall that $f$ is the objective function of problem~\eqref{eq:distributed-problem} and $\delta_{\bar{\mathcal{Y}}}$ is the indicator function of the set $\bar{\mathcal{Y}}$, \cf \cite[p. 6]{Rockafellar-98}.}, where $\bar{\mathcal{Y}}$ is the $m$-fold Cartesian product of $\mathcal{Y}$, i.e.,
\begin{equation}\label{eq:Cartisian-Profuct-Y}
    \bar{\mathcal{Y}}={\mathcal{Y}}^m.
\end{equation} 
\begin{lemma}\label{Lemma:dual-conjugate-relation}
Let $f^*:\R^{nm}\rightarrow \overline{\R}$ denote the conjugate function of $f+\delta_{\bar{\mathcal{Y}}}$. Then $g$ is a restriction of $f^*$ to a linear space. In particular, 
\begin{equation} \label{eq:dual-conjugate-relation}
    g(\boldsymbol{\lambda})=-f^{*}(\mathbf{A}\tran\boldsymbol{\lambda}),
\end{equation}
where $\mathbf{A}$ is an ${n(m-1)\times nm}$ matrix with the special block structure given by: 
\begin{equation}\label{eq:A-Matrix-for-Consensus-Constraint}
\mathbf{A}=\begin{bmatrix} 
\mathbf{I}_n & -\mathbf{I}_{n} & \mathbf{0} & \cdots & \cdots &  \cdots & \mathbf{0}\\  
\mathbf{0} & \mathbf{I}_n & -\mathbf{I}_{n} & \mathbf{0} &  \cdots & \cdots & \mathbf{0} \\ %
\mathbf{0} & \mathbf{0} & \mathbf{I}_n & -\mathbf{I}_{n} &  \mathbf{0} & \cdots & \mathbf{0}\\ %
\vdots & \vdots & \ddots & \ddots & \ddots & \ddots & \vdots   \\  
\mathbf{0} & \vdots & \ddots & \mathbf{0}  &  \mathbf{I}_n & -\mathbf{I}_n & \mathbf{0}\\
\mathbf{0} & \mathbf{0} & \cdots & \cdots  &  \mathbf{0} & \mathbf{I}_n & -\mathbf{I}_n
\end{bmatrix}  .
\end{equation} 
\end{lemma}
%
\begin{IEEEproof}
This is immediate from \eqref{eq:dual-function} and the definition of the conjugate function.
\end{IEEEproof}


\subsection{Lipschitzian Properties}\label{subsecsec:Characteristics-of-Dual-Function-Lipshitzian-Properties}
Let us start by furnishing a simple, yet important result that verifies the Lipschitzian properties of the dual objective function~$g$ of problem~\eqref{eq:distributed-problem}. A hypothesis that is considered in this regard is Assumption~\ref{Assumption:Strongly-Convex-f}, leading to the following lemma.
\begin{lemma}\label{Lemma:Strongly-Convex-f}
Suppose Assumption~\ref{Assumption:Strongly-Convex-f} holds. Then, the objective function $f$ of problem~\eqref{eq:distributed-problem}, given by $f(\mathbf{y})=\sum_{i{=}1}^{m}f_i(\mathbf{y}_i)$, is strongly convex with constant $\mu=\min_i~\mu_i$. 
\end{lemma}
\begin{IEEEproof}
The condition that $f_i$ is strongly convex with constant $\mu_i>0$ is equivalent to the strong monotonicity condition [\cf \cite[Exercise~12.59]{Rockafellar-98}]. The strong convexity of $f$ follows readily by combining the conditions for each $f_i$. Details are omitted due to space limitations. 
\end{IEEEproof}
An immediate consequence of Lemma~\ref{Lemma:Strongly-Convex-f} is the strong convexity of $f+\delta_{\bar{\mathcal{Y}}}$.
The result is outlined in the following remark.
\begin{remark}\label{Remark:Strongly-Convex-f-plus-Delta}
Suppose Assumption~\ref{Assumption:Strongly-Convex-f} holds. Then, the function $f+\delta_{\bar{\mathcal{Y}}}$ is strongly convex with constant $\mu=\min_i~\mu_i$.
\end{remark}
\begin{IEEEproof}
This is immediate from Lemma~\ref{Lemma:Strongly-Convex-f} and the convexity of $\bar{\mathcal{Y}}$, \cf \eqref{eq:Cartisian-Profuct-Y}.
\end{IEEEproof}

The next result claims, under certain assumptions, a useful Lipschitzian property of the dual function~$g$.




\begin{prop}\label{Proposition:Lipscitz-Continuity-of-Grad-g}
Suppose Assumption~\ref{Assumption:Strongly-Convex-f} holds. Then the dual function $g$ is differentiable. Moreover, gradient $\nabla g$ of $g$ is Lipschitz continuous with constant $(1/\mu)\left[2+2\cos(\pi/m)\right]$, where $\mu=\min_i~\mu_i$. 
\end{prop}
\begin{IEEEproof}
Here the closedness of $\mathcal{Y}$ and $f_i$s is vital. The differentiability is a consequence of \cite[Theorem 11.13]{Rockafellar-98}. The Lipschitz continuity is immediate from\cite[Proposition~12.60: (a),(b)]{Rockafellar-98}, together with the biconjugate property $(f^{*})^{*}=f+\delta_{\bar{\mathcal{Y}}}$. The constant is simply the $\lambda_{\max}(\mathbf{A}\mathbf{A}\tran)$, the maximum eigenvalue of $\mathbf{A}\mathbf{A}\tran$, \cf \eqref{eq:A-Matrix-for-Consensus-Constraint} and \cite[p.~565]{Meyer-2000}. Details are omitted due to space limitations.
\end{IEEEproof}

\subsection{Strong Convexity Properties}\label{subsecsec:Characteristics-of-Dual-Function-Strong-Convexity-Properties}
The closedness of $f$, together with the Legendre-Fenchel transform \cite[Theorem 11.1]{Rockafellar-98}, allows a dual result of Proposition~\ref{Proposition:Lipscitz-Continuity-of-Grad-g} to be worked out, again by using \cite[Proposition~12.60: (a),(b)]{Rockafellar-98}. The corresponding results are outlined in the sequel, under Assumption~\ref{Assumption:Lipschitz-Convex-Grad-f}. 
\begin{lemma}\label{Lemma:Lipshitz-Coninnuity-f}
Suppose Assumption~\ref{Assumption:Lipschitz-Convex-Grad-f} holds. Then, the gradient $\nabla f$ of $f$ is Lipschitz continuous on $\R^{nm}$ with constant $L=\max_i~L_i$.  
\end{lemma}

\begin{IEEEproof}
This is clear from straightforward algebraic manipulations and is omitted.
\end{IEEEproof}

The following lemma asserts strongly convexity properties of $-g$ under mild conditions. 

\begin{prop}\label{Proposition:Strong-Concavity-of-g}
Suppose Assumption~\ref{Assumption:Lipschitz-Convex-Grad-f} holds. Then the function $-g$ is strongly convex with constant $(1/L)\left[2-2\cos(\pi/m)\right]$, where $L=\max_i~L_i$. 
\end{prop}
\begin{IEEEproof}
The result follows from Lemma~\ref{Lemma:Lipshitz-Coninnuity-f} and  \cite[Proposition~12.60: (a),(b)]{Rockafellar-98}, together with Lemma~\ref{Lemma:dual-conjugate-relation}. Here again, the closedness of $\mathcal{Y}$ and $f_i$s has a vital role.  Details are omitted due to space restrictions.
\end{IEEEproof}

\subsection{Bounding Properties for the Primal Error}\label{subsecsec:Characteristics-of-Dual-Function-Bounding-Primal-Error-Properties}

Let us finally invoke the strong duality assumption, i.e., Assumption~\ref{Assumption:Strong-Duality}. The following assertions provide a bridge to the quantification of convergence in the primal-domain through related dual results.

\begin{lemma}\label{Lemma:Strong-Duality}
Suppose Assumption~\ref{Assumption:Strongly-Convex-f} and Assumption~\ref{Assumption:Strong-Duality} hold. If the functions $f_i$, $i=1,\ldots,m$ are differentiable, then
\begin{enumerate}
    \item   $ g(\boldsymbol{\lambda}^\star)- g(\boldsymbol{\lambda})\geq {\mu}/{2} \ \|\mathbf{y}(\boldsymbol{\lambda})-\mathbf{y}^\star\|^2$ for all $\boldsymbol{\lambda}\in\R^{n(m-1)}$, 
\item $g(\boldsymbol{\lambda}^\star)-g(\boldsymbol{\lambda})+ S \|\boldsymbol{\lambda}\| \sqrt{g(\boldsymbol{\lambda}^\star)-g(\boldsymbol{\lambda})} \geq \|f(\mathbf{y}(\boldsymbol{\lambda}))-f(\mathbf{y}^\star)\|$, 
\end{enumerate}
where 
\begin{equation}
    \mathbf{y}(\boldsymbol{\lambda}) = \underset{\mathbf{y}\in\bar{\mathcal{ Y}}}{\arg\min}~\left(f(\mathbf{y})+ \boldsymbol{\lambda}\tran \mathbf{A} \mathbf{y}\right),
\end{equation}
the set $\mathcal{\bar Y}$ is given in \eqref{eq:Cartisian-Profuct-Y}, the matrix $\mathbf{A}$ is given in \eqref{eq:A-Matrix-for-Consensus-Constraint}, $\mu=\min_i~\mu_i$, and $S=\sqrt{(4+4\cos(\pi/m))/\mu}$.
\end{lemma}
\begin{IEEEproof}
Let us first define compactly the partial Lagrangian $L:\R^{nm}\times \R^{n(m-1)}\rightarrow\R$ associated with problem~\eqref{eq:distributed-problem} [\cf \eqref{eq:dual-function}], i.e., $L(\mathbf{y},\boldsymbol{\lambda})=f(\mathbf{y})+\boldsymbol{\lambda}\tran \mathbf{Ay}$.
Then,
\begin{align} \label{eq:Bounding-Primal-Error-}
g(\boldsymbol{\lambda}^\star)-g(\boldsymbol{\lambda}) &=  \underset{\mathbf{y}\in\bar{\mathcal{Y}}}{\inf} \ L(\mathbf{y},\boldsymbol{\lambda}^\star)-\underset{\mathbf{y}\in\bar{\mathcal{Y}}}{\inf} \ L(\mathbf{y},\boldsymbol{\lambda}) \\ \label{eq:Bounding-Primal-Error-2}
&= L(\mathbf{y}^\star,\boldsymbol{\lambda}^\star)-L(\mathbf{y}(\boldsymbol{\lambda}),\boldsymbol{\lambda})\allowdisplaybreaks\\ \label{eq:Bounding-Primal-Error-3}
&= L(\mathbf{y}^\star,\boldsymbol{\lambda})-L(\mathbf{y}(\boldsymbol{\lambda}),\boldsymbol{\lambda})\allowdisplaybreaks\\ \label{eq:Bounding-Primal-Error-4}
&\geq \frac{\mu}{2} \ \|\mathbf{y}(\boldsymbol{\lambda})-\mathbf{y}^\star\|^2,\allowdisplaybreaks
\end{align}
where \eqref{eq:Bounding-Primal-Error-} follows from the definition of the dual function, \eqref{eq:Bounding-Primal-Error-2} follows from Assumption~\ref{Assumption:Strong-Duality}, and \eqref{eq:Bounding-Primal-Error-3} is immediate from that $\mathbf{A}\mathbf{y}^{\star}=\mathbf{0}$. Finally, the inequality \eqref{eq:Bounding-Primal-Error-4} follows from \cite[p.~11, (35)]{Polyak-Intro-Opt-1987} since $L$ is a strongly convex function of $\mathbf{y}$ with constant $\mu$, for fixed $\boldsymbol{\lambda}$, and the supposition that $f$ is differentiable. This concludes the proof of the first part. 
For the second claim, we start by \eqref{eq:Bounding-Primal-Error-2}, where we have
\begin{align} \nonumber 
\|g(\boldsymbol{\lambda}^\star)-g(\boldsymbol{\lambda})\| 
&= f(\mathbf{y}^\star) {-}  f(\mathbf{y}(\boldsymbol{\lambda})){-}\boldsymbol{\lambda}\tran \mathbf{A y}(\boldsymbol{\lambda}) \allowdisplaybreaks \\ \nonumber
& \geq \|f(\mathbf{y}(\boldsymbol{\lambda})){-}f(\mathbf{y}^\star)\| {-} \|\boldsymbol{\lambda}\| \| \mathbf{A}\| \|\mathbf{y}(\boldsymbol{\lambda})-\mathbf{y}^\star)\|  \allowdisplaybreaks \\
\nonumber
& = \|f(\mathbf{y}(\boldsymbol{\lambda})){-}f(\mathbf{y}^\star)\| \\ \nonumber
& \qquad -  \sqrt{2+2\cos(\pi/m)} \|\boldsymbol{\lambda}\| \|\mathbf{y}(\boldsymbol{\lambda})-\mathbf{y}^\star)\|.  \allowdisplaybreaks
\end{align}
The first equality follows simply from Assumption~\ref{Assumption:Strong-Duality} and the definition of the Lagrangian. The inequality and the last equality follow immediately by applying the triangular and Cauchy–Schwarz inequalities, together with the properties of the spectral norm of matrix $\mathbf{A}$. The result follows by applying \eqref{eq:Bounding-Primal-Error-4} and rearranging the terms.
\end{IEEEproof}

\section{Problem \eqref{eq:main-problem} over a General Communication Graph}\label{Appendix:General-communication-graph}

Consider a connected graph $G(\mathcal{V},\mathcal{L})$, where the set of nodes $\mathcal{V}=\{1,\ldots,m\}$ represents the SSs of problem~\eqref{eq:main-problem} and the set of edges $\mathcal{L}$ represents direct communication links between nodes. The graph is undirected in the sense that direct communication is permitted between SS~$i$ and SS~$j$ if and only if the same is permitted between SS~$j$ and SS~$i$. To model this, we let the set $\mathcal{L}$ contain both directed edges $(i,j)$ and $(j,i)$ having identical characteristics. We let $\mathcal{N}_i=\{j \ | \ (j,i)\in\mathcal{L}\}$, the set of neighbors of SS~$i$, $i\in \mathcal{V}$. The cardinality of $\mathcal{N}_i$ is denoted by $|\mathcal{N}_i|$. Moreover, let $\mathbf{W}\in\R^{m\times m}$ denote the Laplacian of the graph whose $i$th, $j$th element $w_{ij}$ is given~by
\begin{equation*}\label{eq:Laplacian-matrix}
w_{ij} =
    \begin{cases}
    \displaystyle -1 \ \ \, ; & \ \text{$(i,j)\in\mathcal{L}$}\\
         \displaystyle |\mathcal{N}_i| \,\, ; & \ \text{$i=j$}\\
         \displaystyle 0 \ \ \ \ \, ; & \ \text{$\text{otherwise}$}.
    \end{cases}       
\end{equation*}

We note that $\mathbf{W}$ is symmetric and positive semidefinite. Moreover, with straightforward algebraic manipulations, it is easy to see that $\{\mathbf{y}\in\R^{nm} \ | \  \mathbf{y}_i=\mathbf{y}_{i+1}, \ i\in\mathcal{V}\setminus\{m\}\}=\{\mathbf{y}\in\R^{nm} \ | \  \sqrt{\bar{\mathbf{W}}}\mathbf{y}=\mathbf{0}\}$, where $\mathbf{y}=[\mathbf{y}_1\tran\ldots\mathbf{y}_m\tran]\tran$, $\bar{{\mathbf{W}}}=\mathbf{W}\otimes\mathbf{I}_n$, and $\sqrt{\bar{\mathbf{W}}}$ is a symmetric matrix such that $\sqrt{\bar{\mathbf{W}}}\sqrt{\bar{\mathbf{W}}}=\bar{\mathbf{W}}$.

Thus, problem \eqref{eq:distributed-problem} can be equivalently reformulated as
\begin{equation} \label{eq:distributed-problem-general-graph}
\begin{array}{ll}
\underset{\mathbf{y}}{\mbox{minimize}} & f(\mathbf{y})=\sum_{i{=}1}^{m}f_i(\mathbf{y}_{i}) \\
\mbox{subject to} & \mathbf{y}_{i}\in\mathcal{Y},\ i\in\mathcal{V}\\ 
& \addb{\sqrt{\bar{{\mathbf{W}}}}}\mathbf{y}=\mathbf{0}, 
\end{array}
\end{equation}
where the variable is $\mathbf{y}\in\R^{nm}$. Let $\boldsymbol{\lambda}=[\boldsymbol{\lambda}_{1}\tran\ldots\boldsymbol{\lambda}_{m}\tran]\tran\in\R^{nm}$ denote the Lagrange multiplier vector associated with the constraint $\addb{\sqrt{\bar{{\mathbf{W}}}}}\mathbf{y}=\mathbf{0}$. Then, the dual function of problem \eqref{eq:distributed-problem-general-graph} is given by
\begin{equation*} \label{eq:dual-problem-general-graph}
     g_{\mathbf{W}}(\boldsymbol{\lambda})= \underset{\mathbf{y}_i\in\mathcal{Y}, \ i\in\mathcal{V}}{\inf}\left[f(\mathbf{y})+\boldsymbol{\lambda}\tran\addb{\sqrt{\bar{{\mathbf{W}}}}}\mathbf{y}\right].
\end{equation*}
Hence, the subgradient method to maximize the dual function $g_{\mathbf{W}}$ is given by 
\begin{equation}
   \boldsymbol{\lambda}^{(k+1)}=\boldsymbol{\lambda}^{(k)}+\gamma_k\addb{\sqrt{\bar{{\mathbf{W}}}}}\mathbf{y}^{(k)},\label{eq:general-subgrad-method_2}
\end{equation}
where $\mathbf{y}^{(k)}= {\arg\min}_{\mathbf{y}_i\in\mathcal{Y}, i\in\mathcal{V}}~\left[ f(\mathbf{y})+(\boldsymbol{\lambda}^{(k)})\tran\addb{\sqrt{\bar{{\mathbf{W}}}}}\mathbf{y}\right]$
and $k$ denotes the iteration index. 

\addb{In view of \eqref{eq:general-subgrad-method_2}, the subgradient algorithm does not directly admit a decentralized implementation among SSs. This is because the Lagrangian is not separable, and therefore we cannot evaluate $\mathbf{y}^{(k)}$ over $\mathbf{y}_1,\ldots,\mathbf{y}_m$ separately. However, by introducing $\boldsymbol{\mu}=\sqrt{\bar{{\mathbf{W}}}}\boldsymbol{\lambda}$, where $\boldsymbol{\mu}=[\boldsymbol{\mu}_1\tran\ldots\boldsymbol{\mu}_m\tran]\tran\in\R^{nm}$, a related algorithm to \eqref{eq:general-subgrad-method_2} is obtained, i.e.,
\begin{equation}
   \boldsymbol{\mu}_i^{(k+1)}=\boldsymbol{\mu}_i^{(k)}+\gamma_k\sum_{j\in\mathcal{N}_i\cup \{i\}}w_{ij}\mathbf{y}_j^{(k)}, \ \ i\in\mathcal{V}\label{eq:general-subgrad-method_distributed_2}
\end{equation}
where $\mathbf{y}_j^{(k)}={\arg\min}_{\mathbf{y}_j\in\mathcal{Y}}[ f(\mathbf{y}_j)+(\boldsymbol{\mu}_j^{(k)})\tran\mathbf{y}_j]$ \cite{Sindri-2020}. The update \eqref{eq:general-subgrad-method_distributed_2} turns out to be implementable in a decentralized manner, where the $i$th node needs to communicate only with its neighbors to perform the $\boldsymbol{\mu}_i$~update at each iteration.

It is worth noting that, there is an equivalence between \eqref{eq:general-subgrad-method_2} and \eqref{eq:general-subgrad-method_distributed_2}, if $\boldsymbol{\lambda}^{(0)}$ is chosen from $R(\sqrt{\bar{\mathbf{W}}})$, the range of $\sqrt{\bar{\mathbf{W}}}$ and $\boldsymbol{\mu}^{(0)}$ is chosen from $R({\bar{\mathbf{W}}})$. Moreover, it can be verified that this is not a restriction from the optimality standpoint of the dual problem of \eqref{eq:distributed-problem-general-graph}. Therefore, to maximize the dual function $g_{\mathbf{W}}$, instead of \eqref{eq:general-subgrad-method_2}, one can use \eqref{eq:general-subgrad-method_distributed_2}, which is a decentralized solution method.


}

From~\eqref{eq:Distortion}, together with Assumption~\ref{Assumption:Absolute-Deterministic-Error}, we have $ \|\hat{\mathbf{s}}^{(k)}-\mathbf{s}^{(k)}\|\leq \epsilon_{\mathbf{W}}$, where $\mathbf{s}^{(k)}=\addb{\sqrt{\bar{\mathbf{W}}}}\mathbf{y}^{(k)}$ is the gradient of $g_{\mathbf{W}}$ at $\boldsymbol{\lambda}^{(k)}$ and $\hat{\mathbf{s}}^{(k)}$ is the distorted vector of $\mathbf{s}^{(k)}$. Moreover, we have $\epsilon_{\mathbf{W}}=\addb{\sqrt{\lambda_{\max}(\mathbf{W})\sum_{i=1}^m \varepsilon_i^2}}$, where $\lambda_{\max}(\mathbf{W})$ is the maximum eigenvalue of $\mathbf{W}$.
One can further verify under Assumption~\ref{Assumption:Strongly-Convex-f} that $\nabla g_{\mathbf{W}}$ is Lipschitz continuous with the constant $L_{g_{\mathbf{W}}}=\addb{\lambda_{\max}(\mathbf{W})}/\mu$ [\cf Proposition~\ref{Proposition:Lipscitz-Continuity-of-Grad-g}] and under Assumption~\ref{Assumption:Lipschitz-Convex-Grad-f} that $g_{\mathbf{W}}$ is \addb{strongly concave on $R(\sqrt{\bar{\mathbf{W}}})$} with the constant $\mu_{g_{\mathbf{W}}}=\addb{\lambda_{\min}^{+}(\mathbf{W})}/L$ [\cf Proposition~\ref{Proposition:Strong-Concavity-of-g}], where $\lambda_{\min}^{+}(\mathbf{W})$ is the smallest nonzero eigenvalue of $\mathbf{W}$. Consequently, the theoretical assertions similar to those presented in \S~\ref{sec:Convergence-Analysis-I} and Appendix~\ref{sec:Characteristics-of-Dual-Function} can be derived analogously even when problem~\eqref{eq:main-problem} is considered over a general communication~graph.

\section{Proof of Proposition~\ref{Proposition:Lipschitz-continuous-g-Primal-Result}} \label{Appendix:Proposition:Lipschitz-continuous-g-Primal-Result}

The undistorted local version of the public variable $\mathbf{z}$ is $\mathbf{y}_i\in\R^n$, $i=1,\ldots,m$. Moreover, SSs solve in parallel
\begin{equation} \label{eq:Compact-SubProblem-All-in-One}
\begin{array}{ll}
\underset{\mathbf{y}\in\bar{\mathcal{Y}}}{\mbox{minimize}} & f(\mathbf{y})+ (\boldsymbol{\lambda}^{(k)})\tran \mathbf{A y} 
\end{array}
\end{equation}
to yield the solution $\mathbf{y}^{(k)}=[{\mathbf{y}_1^{(k)\textrm{T}}} \ \ldots \ {\mathbf{y}_m^{(k)\textrm{T}}}]\tran$, \cf line 3 of Algorithm~\ref{Alg:Fully-Distributed}. Then
\begin{align}\label{eq:Lipschitz-continuous-g-Primal-Result}
   \|\mathbf{y}^{(k)}-\mathbf{y}^\star\|^2 &\leq  ({2}/{\mu}) \  \left(h(\boldsymbol{\lambda}^{(k)})-h(\boldsymbol{\lambda}^\star)\right)\\  \allowdisplaybreaks \label{eq:Lipschitz-continuous-g-Primal-Result-1}
   &\leq ({2}/{\mu}) \  \nabla h(\boldsymbol{\lambda}^{(k)})\tran(\boldsymbol{\lambda}^{(k)}-\boldsymbol{\lambda}^\star)\\  \allowdisplaybreaks \label{eq:Lipschitz-continuous-g-Primal-Result-2}
    &\leq ({2}/{\mu}) \  \|\nabla h(\boldsymbol{\lambda}^{(k)})\|\|\boldsymbol{\lambda}^{(k)}-\boldsymbol{\lambda}^\star\|\\  \allowdisplaybreaks \label{eq:Lipschitz-continuous-g-Primal-Result-3}
    & \leq ({2D}/{\mu}) \  \|\nabla h(\boldsymbol{\lambda}^{(k)})\|,
\end{align}
where \eqref{eq:Lipschitz-continuous-g-Primal-Result} follows from the part 1 of Lemma~\ref{Lemma:Strong-Duality}, \eqref{eq:Lipschitz-continuous-g-Primal-Result-1} follows immediately due to the convexity and differentiability of $h$, \eqref{eq:Lipschitz-continuous-g-Primal-Result-2} follows from Cauchy–Schwarz inequality, and finally \eqref{eq:Lipschitz-continuous-g-Primal-Result-3} is directly from the supposition of the proposition. Thus, the part 1 of Corollary~\ref{Corollary:Lipschitz-continuous-nonsum-step} together with \eqref{eq:Lipschitz-continuous-g-Primal-Result-3} yields the first claim. 

Similarly, by using the part 2 of Lemma~\ref{Lemma:Strong-Duality}, we have
\begin{align}\nonumber
  & f(\mathbf{y}^{(k)})-f(\mathbf{y}^\star) \leq  \allowdisplaybreaks \\ \nonumber
  & \hspace{4.5mm}\big[h(\boldsymbol{\lambda}^{(k)}-h(\boldsymbol{\lambda}^\star)\big] + S\|\boldsymbol{\lambda}^{(k)}\|\sqrt{h(\boldsymbol{\lambda}^{(k)})-h(\boldsymbol{\lambda}^\star) }  \allowdisplaybreaks\\ \nonumber
 &\leq    D \|\nabla h(\boldsymbol{\lambda}^{(k)})\| + \sqrt{D}S\|\boldsymbol{\lambda}^{(k)}\| \sqrt{\|\nabla h(\boldsymbol{\lambda}^{(k)})\|}  \allowdisplaybreaks\\ \label{eq:Lipschitz-continuous-g-Primal-Result-second-2}
 &\leq    D \|\nabla h(\boldsymbol{\lambda}^{(k)})\| + \sqrt{D}S(D+\|\boldsymbol{\lambda}^{\star}\|) \sqrt{\|\nabla h(\boldsymbol{\lambda}^{(k)})\|}  \allowdisplaybreaks
\end{align}
where $S{=}\sqrt{(4+4\cos(\pi/m))/\mu}$. Thus, the second part of the proposition follows from the part 1 of Corollary~\ref{Corollary:Lipschitz-continuous-nonsum-step} and \eqref{eq:Lipschitz-continuous-g-Primal-Result-second-2}.

The third and the fourth claims follow from the part 2 of Corollary~\ref{Corollary:Lipschitz-continuous-nonsum-step} together with \eqref{eq:Lipschitz-continuous-g-Primal-Result-3} and \eqref{eq:Lipschitz-continuous-g-Primal-Result-second-2}, respectively. \QE

\section{Proof of Lemma~\ref{Lemma:Feasible-Points}}\label{Appendix:Lemma:Feasible-Points}

Let $\mathbf{B}=(1/m)(\boldsymbol{1}_{m\times m}\otimes \mathbf{I}_n)$ for clarity. First part is immediate from the properties of $\|\mathbf{B}\|$, the spectral norm of~$\mathbf{B}$, together with that $\mathbf{By}^\star=\mathbf{y}^\star$ and $\|\mathbf{B}\|=1$, i.e.,
\begin{align} \label{eq:Bounding-Primal-Error}
\|\mathbf{By}-\mathbf{y}^\star\|^2\leq \|\mathbf{B}\|\| \mathbf{y}-\mathbf{y}^\star\|^2 \leq \|\mathbf{y}-\mathbf{y}^\star\|^2.
\end{align}
The second assertion is essentially based on the convexity of $f$ and \eqref{eq:Primal-Subgrad-Uniform-Bound}. In particular, we have 
\begin{align} \label{eq:Bounding-Primal-Feasible-Error-function}
f(\tilde {\mathbf{y}})-f(\mathbf{y}^\star) 
&\leq \|\tilde {\boldsymbol{\nu}}\|\|\tilde {\mathbf{y}}-\mathbf{y}^\star\|, \quad \forall~\tilde{\mathbf{y}}\in\mathcal{Y}_{\texttt{feas}}, \ \forall~\tilde{\boldsymbol{\nu}}\in\partial f(\tilde {\mathbf{y}}) \allowdisplaybreaks \\
\label{eq:Bounding-Primal-Feasible-Error-function-2}
& \leq \tilde D \ \|\tilde {\mathbf{y}}-\mathbf{y}^\star\|, \quad \forall~\tilde{\mathbf{y}}\in\mathcal{Y}_{\texttt{feas}}, \allowdisplaybreaks 
\end{align}
where $\mathcal{Y}_{\texttt{feas}}=\{\mathbf{y}\in\texttt{dom}~f \ | \ \mathbf{y}\in\bar{\mathcal{Y}}, \ \mathbf{Ay}=\mathbf{0} \}$. \QE

\end{appendices}

\bibliographystyle{IEEEtran}
\bibliography{IEEEabrv,references_paper}
\end{document}